\documentclass{amsart}
\usepackage{color}
 \usepackage[utf8]{inputenc}
 \usepackage[T1,T2A]{fontenc}
\usepackage[russian,english]{babel}
\usepackage{amsmath,amsfonts,amssymb,epsf, euscript,mathtext}
 \usepackage{amsxtra,mathrsfs,textcomp,xcolor}
 \usepackage{amsthm,euscript, bm,graphicx}
\usepackage{slashed}
\usepackage{pdfsync}
\theoremstyle{plain}
\numberwithin{equation}{section}
\makeatletter
\def\@settitle{\begin{center}%
  \baselineskip14\p@\relax
    \bfseries
  \@title
  \end{center}%
  }
\makeatother
\makeatletter
\def\maketitle{\par
  \@topnum\z@ 
  \@setcopyright
  \thispagestyle{firstpage}
    \ifx\@empty\shortauthors \let\shortauthors\shorttitle
  \else \andify\shortauthors
  \fi
  \@maketitle@hook
  \begingroup
  \@maketitle
  \toks@\@xp{\shortauthors}\@temptokena\@xp{\shorttitle}%
  \toks4{\def\\{ \ignorespaces}}
  \edef\@tempa{%
    \@nx\markboth{\the\toks4
      \@nx\MakeUppercase{\the\toks@}}{\the\@temptokena}}%
  \@tempa
  \endgroup
  \c@footnote\z@
  \def\do##1{\let##1\relax}%
  \do\maketitle \do\@maketitle \do\title \do\@xtitle \do\@title
  \do\author \do\@xauthor \do\address \do\@xaddress
  \do\email \do\@xemail \do\curraddr \do\@xcurraddr
  \do\commby \do\@commby
  \do\dedicatory \do\@dedicatory \do\thanks \do\thankses
  \do\keywords \do\@keywords \do\subjclass \do\@subjclass
}
\makeatother
\newcommand{\field}[1]{\ensuremath{\mathbb{#1}}}
\newcommand{\CC}{\field{C}}

\newcommand{\HH}{\field{H}}

\newcommand{\RR}{\field{R}}

\newcommand{\ZZ}{\field{Z}}

\newcommand{\complex}[1]{\mathsf{#1}} 
\newcommand{\EEE}{\complex{E}}

\newcommand{\PPP}{\complex{\Pi}}

\DeclareMathOperator{\im}{Im}
 
 \DeclareMathOperator{\re}{Re}
\DeclareMathOperator{\Tr}{Tr}
\newcommand{\del}{\partial}

\newcommand{\vp}{\varphi}

\newcommand{\curly}[1]{\mathscr{#1}}

\newcommand{\cB}{\curly{B}}

\newcommand{\cF}{\curly{F}}

\newcommand{\cH}{\curly{H}}

\newcommand{\cL}{\curly{L}}

\newcommand{\cU}{\curly{U}}

\newcommand{\bk}{\backslash}

\newcommand{\vep}{\varepsilon}

\newcommand{\vk}{\varkappa}

\newcommand{\NN}{\field{N}}
\newcommand{\QQ}{\field{Q}}

\newcommand{\vphi}{\varphi}

\makeatletter
\def\blfootnote{\xdef\@thefnmark{}\@footnotetext}
\makeatother

\usepackage{hyperref}
\hypersetup{colorlinks,citecolor=blue}

\begin{document}

\title{Etudes of the resolvent}
\author{L.A. Takhtajan}
\date{}
\dedicatory
{To the blessed memory of my teachers, Ludwig Dmitrievich Faddeev and Askold Ivanovich Vinogradov}
\maketitle
\begin{abstract}Based on the notion of the resolvent and on the Hilbert identities, this paper presents a number of classical results in the theory of differential operators and some of their applications to the theory of automorphic functions and number theory from a unified point of view. For instance, for the Sturm-Liouville operator there is a derivation of the Gelfand-Levitan trace formula, and for the one-dimensional Schr\"{o}dinger operator a derivation of  Faddeev's formula for the characteristic determinant and the Zakharov-Faddeev trace  identities. Recent  results on the spectral theory of a certain functional-difference operator arising in conformal field theory are then presented. The last section of the survey is devoted to the Laplace operator on a fundamental domain of a Fuchsian group of the first kind on the Lobachevsky plane. An algebraic scheme is given for proving analytic continuation of the integral kernel of the resolvent of the Laplace operator and the Eisenstein-Maass series. In conclusion, there is a discussion of the relation between the values of the Eisenstein-Maass series at Heegner points and Dedekind zeta-functions of imaginary quadratic fields, and it is explained why pseudo-cuspforms for the case of the modular group do not provide any information about the zeros of the Riemann zeta-function. 
\end{abstract}
\tableofcontents

      \section{Introduction}
      
This survey is an extended and revised version of my talk at the meeting of the Moscow Mathematical Society on April 1, 2014. The reader is presented with a collage of classical results in the theory of differential operators, written from a single viewpoint, together with some applications to automorphic functions and number theory. The relationship between these topics reflects the unity of mathematics, and their choice reflects the tastes and interests of the author, influenced by the traditions of the Leningrad-St.~Petersburg mathematical school. The survey is intended for a broad readership, from specialists in operator theory and functional analysis to algebraic geometers and theoretical physicists. To convey to a modern reader the elegance and beauty of the results, the achievements of the Soviet mathematical school, we chose a neoclassical style of presentation.

Let us describe the contents of this survey in more detail. In \S\ref{Self} we recall the notion of the resolvent, which plays a main role in the theory of self-adjoint operators in a Hilbert space. In \S \ref{Res-spec} we present the Hilbert identities for the resolvent of a self-adjoint operator $A$, and in \S \ref{Det-1} we give the definition of the regularized determinant $\det A$. Section \ref{Sturm-L} is devoted to the classical Sturm-Liouville theory. Thus, in \S \ref{Res-S-L} we briefly recall well-known facts, and in \S \ref{det-S-L} we derive the celebrated Gelfand-Levitan trace formula. Section \ref{Schr-one} is devoted to the presentation of main results for the one-dimensional Schr\"{o}dinger operator $H$. In particular, in \S \ref{Res-Schr} we introduce the Jost solutions and recall the formula for the resolvent of the operator $H$, in \S \ref{Det-char} we derive Faddeev's formula for the regularized determinant $\det (H- \lambda I)$, and in \S \ref{FZ} we present the derivation of the Zakharov-Faddeev trace identities.

Section \ref{Fun-Dif} is based on \cite{TF15} and is devoted to the spectral analysis of a certain functional-difference operator, a special pseudodifferential operator $H$ of infinite order which arises in conformal field theory and in the representation theory of quantum groups. Thus, in \S \ref{W-O} we introduce the self-adjoint Weyl operators $U$ and $V$ in the definition of the operator $H = U + U^{-1} + V$, and in \S \ref{free} we consider the `unperturbed' operator $H_ {0} = U + U^ {-1}$. In \S \ref{H-E} we define a solution to the scattering problem and the Jost functions, and we give an explicit formula for the resolvent of the self-adjoint operator $H$ in the Hilbert space $L^{2}(\RR)$ together with the eigenfunction expansion theorem. Finally, in \S \ref{M-C} we describe functional-difference operators for mirror curves.

The last section,  \S \ref{Laplace}, is devoted to the spectral theory of the Laplace operator $A$ on a fundamental domain of a Fuchsian group of the first kind $\Gamma$ on the Lobachevsky plane $\HH$. Of particular interest here is the case of a non-compact fundamental domain, when $A$ has a continuous spectrum (considered in Faddeev's classical paper \cite{Fad67}). Thus,  in \S \ref{Res-hyp} we follow \cite{Fad67} and give an algebraic scheme for proving a fundamental result on meromorphic continuation with respect to the variable $s$ of the integral kernel of the resolvent $(A-s(1-s) I)^{- 1}$ of the operator $ A $ and of the Eisenstein-Maass series $E(z, s)$, to the domain $0 <\re s \leq 1$. Moreover, in the half-plane $\re s\geq \frac{1}{2} $, the resolvent kernel can have only simple poles and only poles on the line $\re s= \frac{1}{2} $ and the interval $ [\frac{1}{2}, 1]$, while the Eisenstein-Maass series can have only poles on the interval $[\frac{1}{2}, 1]$, which immediately gives the eigenfunction expansion theorem for the operator $A$. As explained in \cite{VKF73}, it follows from the celebrated Selberg trace formula that the regularized determinant of the operator $A$ is expressed in terms of the Selberg zeta-function.

The arithmetic case $\Gamma = \mathrm{PSL}(2, \ZZ)$ --- the modular group --- is considered in \S\S \ref{Pseudo}--\ref{Heegner}. Thus, in \S \ref{Pseudo} we discuss as a curiosity the `sensation' at the end of the 1970s  about the connection of the eigenvalues of the Laplace operator with the zeros of the Riemann zeta function and $L$-series by means of the so-called pseudo-cuspforms. Using the Hilbert identity, we explain why pseudo-cuspforms do not provide any information about the location of these zeros.

Finally, \S\ \ref{Heegner} contains a discussion of the relationship between the values of the Eisenstein-Maass series at Heegner points and the Dedekind zeta-functions of imaginary quadratic fields.   
Using the uniform distribution of the Heegner points in the fundamental domain of the modular group --- the Linnik asymptotics --- we naively `prove' the Riemann hypothesis.
Of course, such an application of Linnik's asymptotics is unacceptable, as is confirmed by an analogue of the classical Vinogradov-Gauss formula in the critical strip, obtained in the paper  \cite{VT-2} by A.I. Vinogradov and the author. Nevertheless, attempts to relate
the Laplace operator, pseudo-cuspforms and the Heegner points to zeros of Dedekind zeta-functions of imaginary quadratic fields continue to this day.
Evidence of this can be found in the papers of Zagier \cite{Zagier} and Colin de Verdi\`{e}r \cite{CV1,CV2} in the early 1980's, as well as in the recent studies by E. Bombieri and P. Garrett (see conference talks \cite{Bom, Gar}).

\section{Main definitions} \label{Self}
In response to questions in quantum mechanics, von Neumann developed a theory of unbounded self-adjoint operators in a Hilbert space.
According to the Dirac-von Neumann axioms (se \cite{T}, for instance), it is self-adjoint operators that correspond to quantum observables, and the simplest of them --- the position and the momentum of a particle --- are described by unbounded operators.

For the convenience of the reader, we follow the classical monograph \cite{AG} and briefly recall the standard notation and basic facts from the theory of self-adjoint operators.

\subsection{Self-adjoint operators}  An operator $A$  with dense domain $D(A)$ in a Hilbert space $\cH$ is said to be symmetric if 
$$(Af,g)=(f,Ag)$$
for all elemends $f,g\in D(A)$, where $(~,~)$ is the inner product in $\cH$. The adjoint operator $A^{*}$ to the densely defined operator $A$ is defined as follows: 
$g\in D(A^{*})$
if there is $g^{\ast}\in\cH$ such that
$$(Af,g)=(f,g^{*})$$
for all $f\in D(A)$ and then $g^{*}=A^{*}g$. The operator $A$ is said to be self-adjoint if $A=A^{*}$. Clearly, every self-adjoint operator is symmetric. 
An operator $A$ is said to be closed if its graph $\Gamma(A)$ --- the set of pairs $\{f,Af\}$ for all 
$f\in D(A)$ --- is a closed subset in $\cH\oplus\cH$; $A$ admits a closure if the closure of $\Gamma(A)$ in  $\cH\oplus\cH$ is the graph of an operator, that is, there is an operator $\bar{A}$ such that $\overline{\Gamma(A)}=\Gamma(\bar{A})$.   A closed operator defined on the whole of $\cH$ is bounded.

A symmetric operator $A$ is said to be essentially self-adjoint if its closure $\bar{A}$ is a self-adjoint operator. A typical example is the operator
$A=id/dx$ acting in the Hilbert space $L^{2}(\RR)$ and defined on the linear space $C_{0}^{\infty}(\RR)$ of smooth functions with compact support. Its closure $\bar{A}$ is a self-adjoint
operator with the domain $D(\bar{A})=W_{2}^{1}(\RR)$, the Sobolev space of absolutely continuous square integrable functions with square integrable derivative. 

\subsection{Resolvent and the spectral theorem} \label{Res-spec} Let $A$ be a closed operator. The values $\lambda\in\CC$, for which an operator\footnote{Here $I$ is the identity operator in $\cH$.} 
$$R_{\lambda}(A)=(A-\lambda I)^{-1}$$ 
(resolvent\footnote{Notation $R(\lambda,A)$ is also used.} of an operator $A$) exists and is defined everywhere on $\cH$, are called regular values. The set $\rho(A)\subseteq\CC$ of regular values is open and is called the resolvent set.  The spectrum of an operator $A$ is the complement to the regular set, $\sigma(A)=\CC\setminus\rho(A)$. For a self-adjoint operator $\sigma(A)\subseteq\RR$.

Resolvent of an operator $A$ satisfies the relation
\begin{equation} \label{Hilbert-1}
R_{\lambda}(A)-R_{\mu}(A)=(\lambda-\mu)R_{\lambda}(A)R_{\mu}(A),\quad\text{where}\quad \lambda,\mu\in\rho(A),
\end{equation}
which generalizes the elementary algebraic formula
$$\frac{1}{a-\lambda}-\frac{1}{a-\mu}=\frac{\lambda-\mu}{(a-\lambda)(a-\mu)}$$
and is called the first Hilbert identity. It follows from \eqref{Hilbert-1} that $R_{\lambda}(A)$ is a holomorphic function on 
 $\rho(A)$ with the values in the Banach algebra $\cL(\cH)$ of bounded operators on $\cH$. Let $A$ and $B$ be closed operators with a common domain.
Their resolvents satisfy the so-called Hilbert second identity
\begin{equation}\label{Hilbert-2}
R_{\lambda}(A)-R_{\lambda}(B)=R_{\lambda}(A)(B-A)R_{\lambda}(B),\quad  \lambda\in\rho(A)\cap\rho(B),
\end{equation}
which generalizes the algebraic formula
$$\frac{1}{a-\lambda}-\frac{1}{b-\lambda}=\frac{b-a}{(a-\lambda)(b-\lambda)}.$$

The spectral theorem of von Neumann is a fundamental fact in the theory of self-adjoint operators acting in a Hilbert space. In particular, for each self-adjoint operator
$A$ there is unique projection-valued countably additive measure $\EEE$, defined on the $\sigma$-algebra $\cB$ of Borel subsets of the real line,  
such that $\EEE(\emptyset)=0$, $\EEE(\RR)=I$, 
$$D(A)=\left\{f\in\cH: \int_{-\infty}^{\infty}\lambda^{2}d(\EEE_{\lambda}f,f)<\infty\right\},$$
and for $f\in D(A)$
$$Af=\int_{-\infty}^{\infty}\lambda d\EEE_{\lambda}f,$$
where $\EEE_{\lambda}=\EEE((-\infty,\lambda))$ and the integral is understood as a limit of Riemann-Stieltjes sums in the strong topology on $\cH$.
The relationship between the projection-valued measure  $\EEE$ and the resolvent $R_{\lambda}$ of $A$ is given by the formula
\begin{equation}\label{Stone}
\lim_{\vep\to 0^{+}}\frac{1}{2\pi i}\int_{a}^{b}(R_{\lambda+i\vep}-R_{\lambda-i\vep})d\lambda=\EEE((a,b))+\frac{1}{2}(\EEE(\{a\})+\EEE(\{b\})),
\end{equation}
sometimes called Stone's formula. This formula is an operator version of the classical Sokhotski–Plemelj formula 
$$\frac{1}{\lambda-i0}-\frac{1}{\lambda+i0}=2\pi i\delta(\lambda)$$
in the theory of distributions, and is a basis of the eigenfunction expansion theorem for differential operators.

\subsection{Determinant of an operator} \label{Det-1} Here we briefly recall the notion of the characteristic determinant of a self-adjoint operator (see \cite{GK}). 
In the simplest case, when $K$ is a compact self-adjoint operator with trace (a trace class or nuclear operator),
the Fredholm determinant is given by the simple formula 
\begin{equation}\label{F-det}
\det (I-\lambda K) =\prod_{i}\left(1-\lambda_{i}\lambda\right),
\end{equation}
where $\lambda_{i}$ are the eigenvalues of an operator $K$, and this is an entire function\footnote{Since $\mathrm{Tr}\,|K|=\sum_{i}|\lambda_{i}|<\infty$.}.
If the operator $K$ is invertible, then from \eqref{F-det} we easily obtain
$$\frac{d}{d\lambda}\log\det(I-\lambda K)=-\mathrm{Tr}\,R_{\lambda}(A),\quad\text{where}\quad A=K^{-1}.$$
 
This formula can be generalized to a wider class of operators. Namely, if the resolvent $R_{\lambda} (A)$ is of trace class,
then the characteristic determinant $\det(A-\lambda I) $ of $A$ is determined (up to a multiplicative constant) from the relation
\begin{equation}\label{Res-det}
\frac{d}{d\lambda}\log\det(A-\lambda I)=-\mathrm{Tr}\,R_{\lambda}(A).
\end{equation}
As we shall see below, this formula makes sense if we understand $\mathrm{Tr}$ to be a properly regularized trace of the resolvent.
If, in addition, $R_{\lambda}(A)$ is an integral operator acting in $L^ {2} (X) $ with integral kernel $R _ {\lambda}(x, y) $ that is continuous on $ X \times X $, where the subset $ X\subset \RR $ is bounded, then by a well-known theorem\footnote{See the monograph \cite[\S III.10]{GK}, as well as \cite[\S 30.5, Theorem 12]{Lax}.}
\begin{equation}\label{Tr-main}
\mathrm{Tr}\,R_{\lambda}(A)=\int_{X}R_{\lambda}(x,x)dx.
\end{equation}
In the case $X=\RR$ the following formula holds:\footnote{It is sufficient to approximate $R_{\lambda}(A)$ by the operators $P_{n}R_{\lambda}(A)P_{n}$, where the $P_{n}$ are the orthogonal projections onto $L^{2}(-n,n)\subset L^{2}(\RR)$, and to use Theorem 6.3 in Chapter III in \cite{GK}.}
\begin{equation}\label{Tr-main-1}
\mathrm{Tr}\,R_{\lambda}(A)=\lim_{n\to\infty}\int_{-n}^{n}R_{\lambda}(x,x)dx.
\end{equation}

A more general way of introducing a regularized determinant is based on the notion of the zeta function of an elliptic operator (see the survey \cite{Schwarz}, as well as \cite{T} and references therein). 
For simplicity, we assume that $A$ is an elliptic operator with a purely discrete spectrum consisting of non-negative eigenvalues $\lambda_{n}$ of finite multiplicity  accumulating to infinity. The zeta function of  $A$ is defined by
 $$\zeta_{A}(s)=\sum_{\lambda_{n}>0}\frac{1}{\lambda_{n}^{s}},$$
 where it is assumed that the series converges absolutely for $ \re s> a $ for some $a> 0$.
Under fairly general assumptions (for example, for the Sturm-Liouville operator considered below), $ \zeta_ {A} (s) $ admits an analytic (meromorphic) continuation to a domain containing the half-plane $ \re s \geq 0 $ and is regular for $ s = 0 $. Then the regularized determinant $ \det A $ of  $ A $ is defined by
$$\det A=\exp\left\{-\zeta'_{A}(0)\right\},$$
where the prime indicates the derivative with respect to $s$.
The characteristic determinant $ \det (A- \lambda I) $ is defined in a similar way, and for many examples this definition is
consistent with the formula \eqref{Res-det}.

\section{Sturm-Liouville problem} \label{Sturm-L} 
\subsection{The resolvent and the eigenfunction expansion} \label{Res-S-L} Following the classic monograph \cite{LS}, we consider the simplest problem of finding all non-trivial solutions of the Sturm-Liouville equation
\begin{equation}\label{S-L}
 -y'' +v(x)y=\lambda y,\quad0\leq x\leq\pi,\\
\end{equation}
with zero boundary conditions
$$y(0)=y(\pi)=0,$$
where $v(x)$ is a continuous real-valued function on the interval 
 $[0,\pi]$. 
The differential operator 
$$\mathcal{L}=-\dfrac{d^{2}}{dx^{2}}+v(x)$$
is symmetric\footnote{The symmetric operator $\mathcal{L}$ acting in $L^{2}(0,\pi)$ has defect indices $(2,2)$, and its self-adjoint extensions are described by the Sturm-Liouville boundary conditions.} on the subspace $C_{0}^{\infty}(0,\pi)$ of smooth functions with compact support. Its Friedrichs extension is the self-adjoint Sturm-Liouville operator  $L$ in $L^{2}(0,\pi)$ with the domain 
$$D(L)=\{y\in W^{2}_{2}(0,\pi): y(0)=y(\pi)=0\},$$ 
where $W^{2}_{2}(0,\pi)$ is the Sobolev space of square-integrable functions on $(0,\pi)$ with square-integrable generalized derivatives up to second order.

Let $y_{1}(x,\lambda)$ and $y_{2}(x,\lambda)$ be solutions of equation \eqref{S-L} with the boundary conditions 
$$y_{1}(0,\lambda)=0,\quad y_{1}'(0,\lambda)=1\quad\text{and}\quad y_{2}(\pi,\lambda)=0,\quad y_{2}'(\pi,\lambda)=1,$$  and let
$d(\lambda)=y_{1}(\pi,\lambda)$. The function $d(\lambda)$ is entire of order $1/2$ with simple zeros $\lambda_{n}$ corresponding to the simple eigenvalues of the operator $L$ and
 tending to infinity.  When $v(x)\in C^{1}(0,\pi)$, the following asymptotics holds:
 \begin{equation}\label{l-as}
 \lambda_{n}=n^{2}+c+O\left(\frac{1}{n}\right),\quad\text{where}\quad c=\frac{1}{\pi}\int_{0}^{\pi}v(x)dx.
 \end{equation}

The operator $L$ has a pure discrete spectrum and the corresponding eigenfunction expansion theorem follows from  \eqref{Stone}. 
Specifically, the resolvent $R_{\lambda}=(L-\lambda I)^{-1}$ of $L$ is an operator-valued meromorphic function with simple poles at $\lambda=\lambda_{n}$ and with residues that are projection operators onto the one-dimensional subspaces corresponding to the eigenfunctions. For $\lambda\neq\lambda_{n}$ the resolvent $R_{\lambda}$ is a bounded operator on $L^{2}(0,\pi)$ with integral
kernel
\begin{equation} \label{formula}
R_{\lambda}(x,\xi)=\frac{1}{W(y_{1},y_{2})(\lambda)}(y_{1}(x,\lambda)y_{2}(\xi,\lambda)\theta(\xi-x)+y_{1}(\xi,\lambda)y_{2}(x,\lambda)\theta(x-\xi)),
\end{equation}
where $W(f,g)=f'g-fg'$ is the Wronskian of functions $f$ and $g$, so that the Wronskian of two solutions of equation \eqref{S-L} does not depend on $x$ and
$W(y_{1},y_{2})(x,\lambda)=-d(\lambda)$, and  $\theta(x)$ is the Heaviside function: $\theta(x)=0$ for $x<0$ and $\theta(x)=1$ for $x\geq 0$.

Indeed, the kernel $R_{\lambda}(x,\xi)$ satisfies the equation
\begin{equation}\label{Res-1}
\left(-\frac{\del^{2}}{\del x^{2}}+v(x)-\lambda\right)R_{\lambda}(x,\xi) 
=\delta(x-\xi),\quad 0<x,\xi<\pi
\end{equation}
(where $\delta(x)$ is the Dirac delta function), which follows from \eqref{S-L} and the elementary formula
$$\theta'(x)=\delta(x)$$ 
in the theory of distributions. Using \eqref{Res-1}, we easily to show that the range of 
 $R_{\lambda}$ is $D(L)$ and $(L-\lambda I)R_{\lambda}=I$. 

\subsection{Characteristic determinant and trace identities}\label{det-S-L}
It follows from \eqref{l-as} that the operator $R_{\lambda}$ is of trace class when $\lambda\neq\lambda_{n}$, so that by using the definition of $\det(L-\lambda I)$ in terms of the operator zeta function, it is not difficult to prove the formula \eqref{Res-det} (e.g., \cite[\S5.5.1]{T}). Since the integral kernel $R_{\lambda}(x,\xi)$ of the trace class operator $R_{\lambda}$ is continuous on $[0,\pi]\times [0,\pi]$, we get from \eqref{Tr-main} that
\begin{equation} \label{Tr-0}
\Tr R_{\lambda}=\int_{0}^{\pi}R_{\lambda}(x,x)dx=-\frac{1}{d(\lambda)}\int_{0}^{\pi}y_{1}(x,\lambda)y_{2}(x,\lambda)dx.
\end{equation}

It is easy to compute the integral in \eqref{Tr-0} using the following classical trick \cite{Fad64}. Specifically, we differentiate \eqref{S-L} for $y_{1}(x,\lambda)$ with respect to $\lambda$:
$$-\dot{y}_{1}''+v(x)\dot{y}_{1}=\lambda\dot{y}_{1}+y_{1},$$
where the dot means the $\lambda$-derivative. We multiply this equation by $y_{2}(x,\lambda)$ and subtract the equation for $y_{2}(x,\lambda)$, multiplied by $\dot{y}_{1}(x,\lambda)$.
As the result, we obtain the identity
\begin{equation}\label{W-0}
y_{1}y_{2}=\dot{y}_{1}y''_{2}- \dot{y}''_{1}y_{2}=-W(\dot{y}_{1},y_{2})'.
\end{equation}
Recalling the definition of the solutions $y_{1}$ and $y_{2}$, we get  from this that
$$\int_{0}^{\pi}y_{1}(x,\lambda)y_{2}(x,\lambda)dx=\dot{y}_{1}(\pi,\lambda)=\dot{d}(\lambda),$$
and comparing \eqref{Res-det}  with \eqref{Tr-0}, we obtain
$$\det(L-\lambda I)=Cd(\lambda),$$
with some constant $C$. By computing the asymptotics as  $\lambda\to-\infty$  in this formula it is easy to get that $C=2$ (see. \cite[\S5.5.1]{T}). For example,
$$\det\left(-\frac{d^{2}}{dx^{2}}-\lambda I\right)=2\frac{\sin\pi\sqrt{\lambda}}{\sqrt{\lambda}}.$$
Summarizing, we obtain the following Hadamard product representation for the entire function $d(\lambda)$:
\begin{equation} \label{Had}
d(\lambda)=\frac{\det L}{2}(-\lambda)^{\delta}\prod_{\lambda_{n}\neq 0}\left(1-\frac{\lambda}{\lambda_{n}}\right),
\end{equation}
where $\delta=1$ if $\lambda=0$ is an eigenvalue for $L$, and $\delta=0$ otherwise.

When $v(x)\in C^{2}(0,\pi)$, one can thoroughly investigate the asymptotics of the function $d(\lambda)$ as $\lambda\to-\infty$ both with the help of the differential equation  \eqref{S-L} and with the help of a Hadamard product and the asymptotics  \eqref{l-as} of the eigenvalues with the remainder term $O(n^{-2})$. Specifically, put $\lambda=-k^{2}$, where $k>0$.  The differential equation \eqref{S-L} with respect to $y_{1}$ is equivalent to the Liouville integral equation
\begin{equation*}
 y_{1}(x,\lambda)=\frac{\sinh kx}{k}+\frac{1}{k}\int_{0}^{x}\sinh\{k(x-t)\}v(t)y_{1}(t,\lambda)dt.
 \end{equation*}
Solving it by the method of successive approximations and integrating by parts, we get after simple calculations that as $k\to \infty$
 \begin{equation}\label{d-as}
 d(\lambda)= \frac{e^{\pi k}}{2k}\left\{1+\frac{\pi c}{2k}+\frac{1}{8k^{2}}\left(\pi^{2}c^{2}-2(v(0)+v(\pi))\right)+O\left(\frac{1}{k^{3}}\right)\right\}.
 \end{equation}
 On the other hand, using the Euler formula for the function $\sinh\pi k$, we rewrite the right hand side of \eqref{Had} as\footnote{Here we assume that $\delta=0$ в \eqref{Had}, which can always be achieved by shifting $v(x)$.} 
 \begin{equation*}
 \Phi(\lambda) =\frac{\det L}{2}\frac{\sinh\pi k}{\pi k}\prod_{n=1}^{\infty}\frac{n^{2}}{\lambda_{n}}\prod_{n=1}^{\infty}\frac{k^{2}+\lambda_{n}}{k^{2}+n^{2}}=C_{1}\frac{\sinh\pi k}{\pi k}\vp(k),
 \end{equation*}
 where
 $$C_{1}= \frac{\det L}{2}\prod_{n=1}^{\infty}\frac{n^{2}}{\lambda_{n}}\quad\text{and}\quad \vp(k)=\frac{\sinh\pi k}{\pi k}\prod_{n=1}^{\infty}\left(1+\frac{\lambda_{n}^{2}-n^{2}}{k^{2}+n^{2}}\right).$$
 Put $s_{\lambda}=\sum_{n=1}^{\infty}(\lambda_{n}-n^{2}-c)$. It follows from \eqref{l-as} with the remainder term $O(n^{-2})$ that
\begin{align*}
\sum_{n=1}^{\infty}\frac{\lambda_{n}^{2}-n^{2}}{k^{2}+n^{2}} & =c\sum_{n=1}^{\infty}\frac{1}{k^{2}+n^{2}}+\frac{1}{k^{2}}s_{\lambda}+O\left(\frac{1}{k^{3}}\right)\\
&=\frac{\pi c\coth\pi k}{2k} -\frac{c}{2k^{2}}+\frac{1}{k^{2}}s_{\lambda}+O\left(\frac{1}{k^{3}}\right).
\end{align*}
From this it is now simple to show (see \cite{LS}) that as $k\to\infty$
\begin{equation}\label{Phi-as}
\Phi(\lambda)= \frac{C_{1}e^{\pi k}}{2\pi k}\left\{1+\frac{\pi c}{2k}+\frac{1}{8k^{2}}\left(\pi^{2}c^{2}-4c+8s_{\lambda}\right)+O\left(\frac{1}{k^{3}}\right)\right\}.
\end{equation}
 
 Comparing the coefficients in the asymptotic formulas \eqref{d-as} and \eqref{Phi-as},  we obtain
 $$C_{1}=\pi\quad\text{и}\quad s_{\lambda}-\frac{c}{2}=-\frac{v(0)+v(\pi)}{4}.$$
The first of these formulas gives the  expression
$$ \det A=2\pi\prod_{n=1}^{\infty}\frac{\lambda_{n}}{n^{2}}=2\pi\prod_{n=1}^{\infty}\left(1+\frac{\lambda_{n}-n^{2}}{n^{2}}\right)$$
 for the regularized determinant of the Sturm-Liouville operator, while the second formula, written in the form\footnote{We take the opportunity to note that this formula corrects a typing error in the corresponding formula in  \S2.2 of  \cite{T17}.}
 $$\sum_{n=1}^{\infty}\left(\lambda_{n}-n^{2}-\frac{1}{\pi}\int_{0}^{\pi}v(x)dx\right)=\frac{1}{2\pi}\int_{0}^{\pi}v(x)dx-\frac{v(0)+v(\pi)}{4},$$
is the celebrated Gelfand-Levitan trace formula  \cite{GL} for the regularized trace of the operator $L$!
 In the case where $v(x)\in C^{\infty}(0,\pi)$, formulas for the regularized traces of all positive-integer powers of the operator $L$ were obtained in Dikii's classical paper \cite{Dik} 
(see also the survey \cite{Sad} and references there).

\section{One-dimensional Schr\"{o}dinger operator}\label{Schr-one}
Leaving aside the case of the radial Schr\"{o}dinger equation (see the survey \cite{Fad59}, the monograph \cite{Mar77} and references there), we consider the Schr\"{o}dinger equation on the whole real line
\begin{equation}\label{Schrod-0}
 -y'' +v(x)y=\lambda y,\quad -\infty<x<\infty.
 \end{equation}
 Here the potential --- a measurable, real-valued function $v(x)$ --- is assumed to be rapidly decaying:
\begin{equation} \label{v-small}
\int_{-\infty}^{\infty}(1+|x|)|v(x)|dx<\infty.
\end{equation}
Without loss of generality, we assume $v(x)$ to be continuous. 
Under condition \eqref{v-small}, the Schr\"{o}dinger operator 
$$H=-\dfrac{d^{2}}{dx^{2}}+v(x)$$
is defined on the functions $\psi\in L^{2}(\RR)$ that are twice differentiable on 
$\RR$  and such that $-\psi''+v(x)\psi\in L^{2}(\RR)$, and it is self-adjoint in $L^{2}(\RR)$. 
It is convenient to write
$$H=H_{0}+V,$$
where $H_{0}=-d^{2}/dx^{2}$ is the free Schr\"{o}dinger operator with $D(H_{0})=W^{2}_{2}(\RR)$, and $V$ is a multiplication by $v(x)$ operator in $L^{2}(\RR)$. The operator $H$ has an absolutely continuous spectrum of multiplicity $2$ filling the semi-axis 
$[0,\infty)$, and finitely many simple negative eigenvalues $\lambda_{1},\dots,\lambda_{n}$. Let us explain this more carefully (see  \cite{Fad64,Fad74,Mar77,T} for details).

\subsection{Jost solutions and the resolvent}\label{Res-Schr}
It is convenient to use the parametrization 
$\lambda=k^{2}$, in which the complex $\lambda$-plane cut along  
$[0,\infty)$ corresponds to the upper half-plane of the variable $k$, the so-called `physical sheet' of the Riemann surface of the function 
 $k=\sqrt{\lambda}$. Under the condition \eqref{v-small} the Jost solutions are defined, namely,  
 the functions $f_{1}(x,k)$ and $f_{2}(x,k)$ satisfying \eqref{Schrod-0}  and having the asymptotics 
\begin{alignat*}{3}
f_{1}(x,k) &=e^{ikx}+o(1), &\quad \text{as}\quad &x\to\infty, \\
f_{2}(x,k) &=e^{-ikx}+o(1), &\quad \text{as}\quad & x\to-\infty.
\end{alignat*}
The proof is based on the integral equation of the Volterra type
$$f_{1}(x,k)=e^{ikx}-\int_{x}^{\infty}\frac{\sin k(x-t)}{k}v(t)f_{1}(t,k)dt,$$
and the analogous equation for $f_{2}(x,k)$. The following estimate holds for $k\neq 0$:
\begin{equation}\label{J-0}
|e^{-ikx}f_{1}(x,k)-1|\leq \frac{\sigma(x)}{k}e^{\frac{1}{|k|}\sigma(x)},\quad\text{where}\quad\sigma(x)=\int_{x}^{\infty}|v(t)|dt,
\end{equation}
and moreover $\sigma\in L^{1}(a,\infty)$ for each $a\in\RR$ and $\lim_{x\to\infty}\sigma(x)=0$. 
Similarly,
\begin{equation}\label{J-1}
|e^{ikx}f_{2}(x,k)-1|\leq \frac{\tilde\sigma(x)}{k}e^{\frac{1}{|k|}\tilde\sigma(x)},\quad\text{where}\quad\tilde\sigma(x)=\int_{-\infty}^{x}|v(t)|dt,
\end{equation}
and moreover $\tilde\sigma\in L^{1}(-\infty, a)$ for each $a\in\RR$ and $\lim_{x\to-\infty}\tilde\sigma(x)=0$. 

For fixed $x$ the Jost solutions $f_{1}(x,k)$ and
 $f_{2}(x,k)$ admit analytic continuation to the  half-plane $\im k>0$ and for fixed $k$ satisfy the estimates \eqref{J-0}--\eqref{J-1}.
For real $k$
\begin{equation} \label{a-b}
f_{2}(x,k)=a(k)f_{1}(x,-k)+b(k)f_{1}(x,k),
\end{equation}
where $\overline{a(k)}=a(-k)$, $\overline{b(k)}=b(-k)$ and 
$$|a(k)|^{2}=1+|b(k)|^{2}.$$
The functions $a(k)$ and $b(k)$ are called transition coefficients\footnote{In quantum mechanics the functions $t(k)=\frac{1}{|a(k)|^{2}}$ and $r(k)=\frac{|b(k)|^{2}}{|a(k)|^2}$ are called the transmission and reflection coefficients, respectively.}. For the coefficient $a(k)$ we have the formula
\begin{equation}\label{a}
a(k)=\frac{1}{2ik}W(f_{1}(x,k),f_{2}(x,k)),
\end{equation}
which implies that  $a(k)$ admits analytic continuation to the upper half-plane 
$\im k>0$ and satisfies
$$a(k)=1+O\left(\frac{1}{|k|}\right)\quad\text{as}\quad k\to\infty.$$ 
The function $a(k)$ in the upper half-plane $\im k>0$ has finitely many simple zeros $i\varkappa_{j}$ on the imaginary semi-axis, and 
$\lambda_{j}=-\varkappa_{j}^{2}$ are the eigenvalues of the operator $H$ with the eigenfunctions 
$\psi_{j}(x)=f_{1}(x,i\varkappa_{j})$, $j=1,\dots,n$. Furthermore, it follows from the Poisson-Schwarz formula 
that $a(k)$ satisfies the so-called dispersion relation
\begin{equation}\label{disp}
a(k)=\exp\left\{\frac{1}{\pi i}\int_{-\infty}^{\infty}\frac{\log |a(q)|}{q-k}dq\right\}\prod_{j=1}^{n}\frac{k-i\varkappa_{j}}{k+i\varkappa_{j}},\quad \im k>0.
\end{equation}

The resolvent $R_{\lambda}=(H-\lambda I)^{-1}$ of the Schr\"{o}dinger operator $H$ is defined on 
$$\rho(H)=\CC\setminus\{[0,\infty)\cup\{\lambda_{1},\dots,\lambda_{n}\}\}$$ 
and is an integral operator in $L^{2}(\RR)$ with the integral kernel  
\begin{equation} \label{Res-Scr}
R_{\lambda}(x,y)=-\frac{1}{2ika(k)}(f_{1}(x,k)f_{2}(y,k)\theta(x-y)+f_{1}(y,k)f_{2}(x,k)\theta(y-x)),
\end{equation}
where $k=\sqrt{\lambda}$. 
In particular, the integral kernel of the resolvent $R^{0}_{\lambda}$ of the free operator $H_{0}$ takes the form
\begin{equation}\label{Res-0}
R^{0}_{\lambda}(x,y)=-\frac{e^{ik|x-y|}}{2ik},\quad \im k>0.
\end{equation}
As in the case of the Sturm-Liouville operator,  
the integral kernel $R_{\lambda}(x,y)$ satisfies the same equation \eqref{Res-1},
\begin{equation}\label{Res-2}
\left(-\frac{\del^{2}}{\del x^{2}}+v(x)-\lambda\right)R_{\lambda}(x,y) =\delta(x-y), 
\end{equation}
where now $-\infty<x,y<\infty$.

The eigenfunction expansion for the operator $H$ follows from the formulas \eqref{Stone} and \eqref{a-b}--\eqref{Res-Scr}. In particular, denote by $P$ the orthogonal projection from $\cH=L^{2}(\RR)$ onto the subspace spanned by $\psi_{1},\dots,\psi_{n}$, and denote by  $\mathfrak{H}_{0}$ the Hilbert space $L^{2}\left([0,\infty),\CC^{2};|a(k)|^{-2}dk\right)$. The operator
$\cU: \cH \rightarrow \mathfrak{H}_{0}$  
defined by the formula
$$(\cU f)_{l}(k)=\frac{1}{\sqrt{2\pi}}\int_{-\infty}^{\infty}f(x)f_{l}(x,k)dx,\quad l=1,2,$$
is a partial isometry of the Hilbert spaces
$\cH$ and $\mathfrak{H}_{0}$:
$$\cU^{*}\cU=I-P,\quad \cU\cU^{*}=I_{0},$$
where $I_{0}$ is the identity operator on $\mathfrak{H}_{0}$. The operator $\cU H\cU^{*}$ is the multiplication by $k^{2}$ operator in $\mathfrak{H}_{0}$. The eigenfunction expansion for the free Schr\"{o}dinger operator is the Fourier transform. 
\subsection{The characteristic determinant}\label{Det-char}
Since the operator $H$ has an absolutely continuous spectrum, the formula \eqref{Res-det} no longer makes sense, and there now arises the problem of defining $\det(H- \lambda I)$. 
In the case of the radial Schrödinger operator, this problem was solved by Buslaev and Faddeev in \cite{Bus-Fad}, which subsequently led to the concept of the perturbation determinant \cite {GK}. Here we consider the case of the one-dimensional Schrödinger operator and, for simplicity of presentation, instead of \eqref{v-small} we impose a stronger condition on the potential $v (x)$.

Namely, suppose that $v(x)$ is a bounded function on the real axis and
$v(x)=O(|x|^{-3-\vep})$ as $|x|\to\infty$ 
for some $\vep>0$. It follows from the first condition that $V$ is a bounded operator in $L^{2}(\RR)$, and the second condition means that in the estimates  \eqref{J-0}--\eqref{J-1} one can replace $
\sigma(x)$ and $\tilde\sigma(x)$ by $O(|x|^{-2-\vep})$. By analogy with \eqref{Res-det}, the regularized determinant of the operator $H-\lambda I$ is given by
\begin{equation}\label{det-H}
-\frac{d}{d\lambda}\log\det(H-\lambda I)=\mathrm{Tr}(R_{\lambda}-R_{\lambda}^{0}),\quad \lambda\in\rho(H),
\end{equation}
where $R_{\lambda}-R_{\lambda}^{0}$ is a trace class operator. Indeed,  from the second Hilbert identity we obtain
\begin{equation} \label{Hilb-2}
R_{\lambda}-R_{\lambda}^{0}=R_{\lambda}VR_{\lambda}^{0}.
\end{equation}
Denote by $\sqrt{V}$ the multiplication operator by the function $\sqrt{v(x)}$,
where 
$$\sqrt{v(x)}=\sqrt{|v(x)|}e^{i\delta},$$ 
$\delta=0$ if $v(x)\geq 0$ and $\delta=i\pi/2$ if $v(x)<0$. 
Since $\sqrt{v(x)}\in L^{2}(\RR)$, it follows from \eqref{Res-Scr} and \eqref{Res-0} that the operators $R _ {\lambda} \sqrt{V} $ and $\sqrt{V} R^{0}_{\lambda} $ are Hilbert-Schmidt operators, and therefore the operator $R_{\lambda} -R_{\lambda}^{0} $ is of trace class.
 
Remarkably, the trace on the right-hand side of the formula \eqref{det-H} can be calculated explicitly. In case of the radial Schrödinger equation, the corresponding formula was given by
Buslaev and Faddeev in \cite{Bus-Fad}, and in case of the whole axis, this is the last formula in section \S1.1 of Faddeev's survey \cite{Fad74} (see also Problem 2.6 in \S3.2.2 of chapter 3 in \cite{T}). In particular, the following relation holds:
\begin{equation} \label{det-F}
\mathrm{Tr}(R_{\lambda}-R_{\lambda}^{0})=-\frac{d}{d\lambda}\log a(\sqrt{\lambda}).
\end{equation}
As far as we know, no complete derivation of this beautiful formula exists in the literature. For the convenience of the reader, we present it here.

\begin{proof} Recall that $k=\sqrt{\lambda}$. For $\im k>0$ the Jost solution $f_{1}(x,k)$ decays exponentially for large  $x$. For such $k$ a solution $g(x,k)$ of \eqref{Schrod-0} linear independent from
 $f_{1}(x,k)$ is found from the relation
$W(f_{1},g)=2ik$, and it grows exponentially for large $x$:
$$e^{ikx}g(x,k)=1+O(x^{-1-\vep})\quad\text{and}\quad e^{ikx}g'(x,k)=-ik+O(x^{-1-\vep})\quad\text{as}\quad x\to\infty.$$
The functions $f_{1}(x,k)$ and $g(x,k)$ form a basis in the solution space, and from the condition $W(f_{1},g)=2ik$ and \eqref{a}  we get that
$$f_{2}(x,k)=a(k)g(x,k)+c(k)f_{1}(x,k).$$ Therefore, as $x\to\infty$ we have
\begin{equation}\label{2-grow}
e^{ikx}f_{2}(x,k)=a(k)+O(x^{-1-\vep})\quad\text{and}\quad e^{ikx}f_{2}'(x,k)=-ika(k)+O(x^{-1-\vep}). 
\end{equation}
Similarly, as $x\to-\infty$
\begin{equation}\label{1-grow}
e^{-ikx}\!f_{1}(x,k)=a(k)+O(|x|^{-1-\vep})\;\;\text{and}\;\; e^{-ikx}f_{1}'(x,k)=ika(k)+O(|x|^{-1-\vep}).
\end{equation}

The trace class operator $R_{\lambda}-R^{0}_{\lambda}$ has an  integral kernel  $\tilde{R}_{\lambda}(x,y)$ that is continuous on  $\RR\times\RR$ and
$$\tilde{R}_{\lambda}(x,x)=-\frac{1}{2ika(k)}\left(f_{1}(x,k)f_{2}(x,k)-a(k)\right).$$
Thus, using formula \eqref{Tr-main-1}, we obtain
\begin{equation}\label{tr-dif}
\Tr(R_{\lambda}-R^{0}_{\lambda})=
\frac{i}{2ka(k)}\lim_{n\to\infty}\int_{-n}^{n}\left(f_{1}(x,k)f_{2}(x,k)-a(k)\right)dx.
\end{equation}
As in the case of a Sturm-Liouville operator, the integral in \eqref{tr-dif} can be evaluated explicitly. Namely, write \eqref{W-0} in the form
\begin{equation}\label{W-1}
f_{1}(x,k)f_{2}(x,k)=-\frac{1}{2k}W(\dot{f}_{1}(x,k), f_{2}(x,k))'=\frac{1}{2k}W(f_{1}(x,k),\dot{f}_{2}(x,k))',
\end{equation}
where the dot now stands for the $k$-derivative. From the integral equation for  $f_{1}(x,k)$ we obtain
$$e^{-ikx}\dot{f}_{1}(x,k)=ix+O(x^{-\vep})\quad\text{and}\quad e^{-ikx}\dot{f}'_{1}(x,k)=i-kx+O(x^{-\vep})\quad\text{as}\quad x\to\infty,$$
and therefore for such $x$ we get by using \eqref{2-grow} that
$$W(\dot{f}_{1}(x,k), f_{2}(x,k))=(i-2kx)a(k)+O(x^{-\vep}).$$
Similarly, as $x\to-\infty$,
$$W(f_{1}(x,k), \dot{f}_{2}(x,k))=(i+2kx)a(k)+O(|x|^{-\vep}).$$
Using these formulas and \eqref{W-1}, we obtain
\begin{align*}
\int_{0}^{n}f_{1}(x,k)f_{2}(x,k)dx & =\left(-\frac{i}{2k}+n\right)a(k) +\frac{1}{2k}W(\dot{f}_{1}(0,k),f_{2}(0,k))+ O(n^{-\vep})\\
\intertext{and}
\int_{-n}^{0}f_{1}(x,k)f_{2}(x,k)dx & =\left(-\frac{i}{2k}-n\right)a(k) +\frac{1}{2k}W(f_{1}(0,k),\dot{f}_{2}(0,k))+ O(n^{-\vep}).
\end{align*}
Adding two last formulas and using the relation
$$\dot{a}(k)=-\frac{a(k)}{k}+\frac{1}{2ik}W(\dot{f}_{1}(0,k),f_{2}(0,k))+\frac{1}{2ik}W(f_{1}(0,k),\dot{f}_{2}(0,k))$$
following from  \eqref{a},  we have
$$\lim_{n\to\infty}\int_{-n}^{n}\left(f_{1}(x,k)f_{2}(x,k)-a(k)\right)dx=i\dot{a}(k).$$
Substitution of this relation into \eqref{tr-dif} gives the desired formula \eqref{det-F}.
\end{proof}

\subsection{Trace identities}\label{FZ}
It follows from \eqref{det-F} that the regularized determinant of the Schr\"{o}dinger operator is given by the formula 
$$\det(H-\lambda I)=a(\sqrt{\lambda}),$$
is a holomorphic function on the complex $\lambda$-plane cut along the non-negative semi-axis, and 
has zeros at the eigenvalues of the operator $H$. Under the assumption that the potential $v(x)$ is a function in the Schwartz class\footnote{That is, it is smooth and rapidly decaying with all derivatives as $|x|\to\infty$.}, we easily obtain from  \eqref{det-F}, (as in \S
 \ref{Sturm-L}) the trace identities for the one-dimensional Schr\"{o}dinger operator $H$.

Namely, for such $v(x)$ the coefficient $b(k)$  is a Schwartz class function, and therefore from \eqref{disp} we immediately obtain an asymptotic expansion as $k\to\infty$:
\begin{equation} \label{a-as}
\log a(k)=\sum_{l=1}^{\infty}\frac{c_{l}}{k^{l}}+O(|k|^{-\infty}),\quad \im k>0,
\end{equation}
which is an analogue of the expansion \eqref{Phi-as} for the characteristic determinant of the Sturm-Liouville operator. 
Moreover, $c_{2l}=0$ and
\begin{equation}\label{c-as}
c_{2l+1}=-\frac{1}{\pi i}\int_{-\infty}^{\infty}k^{2l}\log|a(k)|dk-\frac{2}{2l+1}\sum_{j=1}^{n}(i\varkappa_{j})^{2l+1}.
\end{equation}

An analogue of the asymptotic expression \eqref{d-as} is obtained by means of the following beautiful argument (here we follow the famous paper  \cite{ZF} by Zakharov and Faddeev).  
From formula \eqref{J-0} it follows that the function $\chi(x,k)=\log f_{1}(x,k)$ is well-defined for large $k$ with $\im k>0$, and
$$\chi(x,k)=ikx+o(1),\quad x\to\infty\quad\text{and}\quad  \chi(x,k)=\log a(k) +ikx +o(1),\quad x\to-\infty.$$
As follows from \eqref{Schrod-0}, the function 
$$\sigma(x,k)=\frac{d}{dx}\chi(x,k)-ik,$$
is a solution of the Riccati equation
$$\sigma'+\sigma^{2}-v+2ik\sigma=0,$$
decays as $|x|\to\infty$ and satisifes
\begin{equation} \label{log-a}
\log a(k)=-\int_{-\infty}^{\infty}\sigma(x,k)dx.
\end{equation}
Now it is not difficult to verify the asymptotic expansion
\begin{equation}\label{sigma-as}
\sigma(x,k)=\sum_{l=1}^{\infty}\frac{\sigma_{l}(x)}{(2ik)^{l}} +O(|k|^{-\infty}).
\end{equation}
The coefficients $\sigma_{l}(x)$ are polynomials in the function $v(x)$ and its derivatives at $x$, and are determined  
by the recurrence relation 
$$\quad \sigma_{l}(x)=-\sigma'_{l-1}(x)-\sum_{j=1}^{l-1}\sigma_{l-j-1}(x)\sigma_{j}(x),\quad \sigma_{1}(x)=v(x);$$
moreover, the $\sigma_{2l}(x)$ are total $x$-derivatives. Comparing the formulas \eqref{a-as}, \eqref{c-as} with \eqref{log-a}, \eqref{sigma-as}, we obtain the Zakharov-Faddeev trace identities
$$\frac{1}{\pi i}\int_{-\infty}^{\infty}k^{2l}\log|a(k)|dk+\frac{2}{2l+1}\sum_{j=1}^{n}(i\varkappa_{j})^{2l+1}=\left(\frac{1}{2i}\right)^{2l+1}\int_{-\infty}^{\infty}\sigma_{2l+1}(x)dx.$$
In \cite{ZF}, the reader can find a remarkable application of these formulas to the proof of complete integrability of the Korteweg-de Vries equation.

Comparing \eqref{det-F} with \eqref{a-as}, \eqref{c-as}, we see that  $\Tr(R_{\lambda}-R^{0}_{\lambda})$ can be expanded as $\lambda\to-\infty$  in an asymptotic series in inverse odd powers of $\sqrt{\lambda}$. 
In the Gel'fand--Dikii paper \cite{GD} this was proved directly, in both the rapidly decreasing case and the periodic case.
Namely, rewriting the second Hilbert identity \eqref{Hilb-2} in the form
$$R_{\lambda}(I-VR_{\lambda}^{0})=R_{\lambda}^{0},$$
we obtain
$$R_{\lambda}=R_{\lambda}^{0}+\sum_{n=1}^{\infty}R_{\lambda}^{0}(VR_{\lambda}^{0})^{n},$$
where the infinite series is understood as an asymptotic series as $\lambda\to-\infty$. By using the explicit formula \eqref{Res-0} 
for the free resolvent, it is not difficult to obtain the asymptotic expansion 
$$R_{\lambda}(x,x)=\sum_{l=1}^{\infty}\frac{R_{l}(x)}{\lambda^{l+\frac{1}{2}}}+O(|\sqrt\lambda|^{-\infty}).$$
The coefficients $R_{l}(x)$ are easily found from the third-order differential equation
$$\left(-\frac{d^{3}}{dx^{3}} +4(v(x)-\lambda)\frac{d}{dx}+2v'(x)\right)R_{\lambda}(x,x)=0,$$
for the product of two solutions of the second-order equation \eqref{Schrod-0}. Details of these beautiful calculations can be found in \cite{GD}.

\section{A certain functional-difference operator} \label{Fun-Dif}
Consider the following functional-difference equation
\begin{equation*}
\psi(x+ib)+\psi(x-ib) +e^{2\pi bx}\psi(x)=\lambda\psi(x), 
\end{equation*}
where 
$$b>0\;\;\text{and}\;\;-\infty<x<\infty,$$
and the function $\psi(x)$ admits analytic continuation into the strip 
$$\Pi_{b}=\{z=x+iy\in\CC: |y|<b\}.$$ 
A functional-difference operator
$$H=U+U^{-1}+V$$
is associated with this equation, 
where $U$ and $V$ are the self-adjoint Weyl operators acting in $L^{2}(\RR)$. The operator $H$ arises in conformal field theory and in the representation theory of the quantum group $\mathrm{SL}_{q}(2,\RR)$.  In \cite{TF15} there is a spectral analysis of this unbounded self-adjoint operator acting in $L^{2}(\RR)$. We give a detailed presentation of these results.
\subsection{Weyl operators} \label{W-O}
The quantum mechanical Weyl operators are unitary operators $U(u)$ and $V(v)$ in $L^{2}(\RR)$,   $u,v\in\RR$, defined by the formulas 
$$(U(u)\psi)(x)=\psi(x-u),\quad (V(v)\psi)(x)=e^{-ivx}\psi(x),\quad \psi\in L^{2}(\RR)$$
(see, e.g., \cite[Ch. 2]{T}, where the Planck constant $\hbar$ is set to be $1$).
The operators $U(u)$ and $V(v)$ satisfy the Weyl commutation relations 
$$U(u)V(v)=e^{iuv}V(v)U(u).$$
In the representation theory of the quantum group $\mathrm{SL}_{q}(2,\RR)$ one uses complex values of $u$ and $v$, under which the Weyl operators $U(u)$ and $V(v)$ become unbounded self-adjoint operators acting in $L^{2}(\RR)$. 

Namely, consider the operators $U$ and $V$, given formally by 
\begin{equation} \label{U-V}
(U\psi)(x)=\psi(x+ib), \quad (V\psi)(x)=e^{2\pi bx}\psi(x)
\end{equation}and satisfying the relation 
\begin{equation} \label{U-V-q}
UV=q^{2}VU,\quad q=e^{\pi ib^{2}}
\end{equation}
on the common domain of $U$ and $V$. 
The operators $U$ and $V$, defined by \eqref{U-V} are unbounded self-adjoint operators acting in $L^{2}(\RR)$. Specifically, $U$  is a self-adjoint operator acting in $L^{2}(\RR)$ with the domain
$$D(U)=\{\psi(x)\in L^{2}(\RR) : e^{-2\pi b p}\hat{\psi}(p)\in L^{2}(\RR)\},$$ 
where
$$\hat{\psi}(p)=\cF(\psi)(p)=\int_{-\infty}^{\infty}\psi(x)e^{-2\pi ipx}dx$$
is the Fourier transform\footnote{We are using the normalization of the Fourier transform that is customary in analytic number theory.} in $L^{2}(\RR)$.  Equivalently, the domain $D(U)$ consists of the functions $\psi(x)$ which admit analytic continuation into the strip 
$$\Pi^{+}_{b}=\{z=x+iy\in\CC: 0< y <b\}$$ 
with the property $\psi(x+iy)\in L^{2}(\RR)$ for all $0\leq y<b$ and that the limit 
$$\psi(x+ib-i0)=\lim_{\varepsilon\rightarrow 0^{+}}\psi(x+ib-i\varepsilon)$$
exists in the sense of convergence in $L^{2}(\RR)$. Furthermore, for $\psi\in D(U)$ we have
$$(U\psi)(x)=\psi(x+ib-i0).$$ 

The domain $D(U^{-1})$ of the inverse operator $U^{-1}$ is defined similarly, and  we have $(U^{-1}\psi)(x)=\psi(x-ib+i0)$. The domain $D(V)$ of the self-adjoint operator $V$ consists of the functions $\psi(x)\in L^{2}(\RR)$ for which $e^{2\pi bx}\psi(x)\in L^{2}(\RR)$. 
Thus, we have
$$U^{-1}=\cF^{-1}V\cF,$$
where the inverse Fourier transform is given by the formula
$$\psi(x)=\int_{-\infty}^{\infty}\hat{\psi}(p)e^{2\pi ipx}dp.$$

\subsection{The operator $H_{0}$} \label{free}
The free operator $H_{0}=U+U^{-1}$ is an unbounded self-adjoint operator acting in $L^{2}(\RR)$, and defined on $D(H_{0})=D(U)\cap D(U^{-1})$ by the formula 
$$(H_{0}\psi)(x)=\psi(x+ib-i0)+\psi(x-ib+i0),\quad \psi\in D(H_{0}).$$
Obviously, for $b\rightarrow 0$ the operator $b^{-2}(H_{0}-2I)$ turns into the operator $-d^{2}/dx^{2}$.
In terms of the Fourier transform the operator $\hat{H}_{0}=\cF H_{0}\cF^{-1}$ is the multiplication  by $2\cosh(2\pi bp)$ operator, and thus domain $D(H_{0})$ admits an equivalent description: 
$$D(H_{0})=\left\{\psi(x)\in L^{2}(\RR) : \int_{-\infty}^{\infty}\cosh^{2}(2\pi bp)|\hat{\psi}(p)|^{2}dp<\infty\right\},$$
and it is a `hyperbolic analogue' of the Sobolev space $W^{2}_{2}(\RR)$. 

For$\lambda\in\CC\setminus [2,\infty)$ the resolvent of the operator $\hat{H}_{0}$,
$$\hat{R}_{\lambda}^{0}=(\hat{H}_{0}-\lambda I)^{-1},$$
is the multiplication by $(2\cosh(2\pi bp)-\lambda)^{-1}$ operator and it is bounded on $L^{2}(\RR)$. Because the function $2\cosh(2\pi bp)$ is a two-to-one map of the real axis $-\infty<p<\infty$ onto $[2,\infty)$, 
the spectrum of $\hat{H}_{0}$ is absolutely continuous and fills the semi-infinite interval $[2,\infty)$ with multiplicity  $2$.
Correspondingly, for $\lambda\in\CC\setminus [2,\infty)$ the resolvent  
$$R_{\lambda}^{0}=(H_{0}-\lambda I)^{-1}$$ 
of $H_{0}$  is an integral operator acting in $L^{2}(\RR)$ with integral kernel depending on the difference of the arguments,
\begin{equation} \label{R-0-def}
(R_{\lambda}^{0}\psi)(x)=\int_{-\infty}^{\infty}R^{0}_{\lambda}(x-y)\psi(y)dy,
\end{equation}
where
\begin{equation} \label{R-0-kernel}
R^{0}_{\lambda}(x)=\int_{-\infty}^{\infty}\frac{e^{2\pi ipx}}{2\cosh(2\pi bp)-\lambda}dp.
\end{equation}
It is convenient to use the parametrization (cf. \S \ref{Res-Schr})
$$\lambda=2\cosh(2\pi bk),$$
in which the resolvent set $\CC\setminus [2,\infty)$ becomes the `physical sheet' --- the strip $0<\im k\leq 1/(2b)$ ---
and the continuous spectrum $[2,\infty)$ is doubly covered by the real axis $-\infty<k<\infty$.  The integral \eqref{R-0-kernel} 
is easily calculated, and we obtain 
\begin{equation} \label{R-0-formula}
R^{0}_{\lambda}(x) =\frac{i}{2b\sinh(2\pi bk)}\left(\frac{e^{-2\pi ikx}}{1-e^{2\pi x/b}}+\frac{e^{2\pi ikx}}{1-e^{-2\pi x/b}} \right),\quad \lambda=2\cosh(2\pi bk).
\end{equation}
The function $R^{0}_{\lambda}(x)$ is regular at $x=0$ and for $0<\im k\leq 1/(2b)$ the following estimate holds:
$$|R^{0}_{\lambda}(x)|\leq Ce^{-2\pi\im k|x|},$$
where $C>0$ is a constant\footnote{Here and below  we use $C$ to denote various constants.}, 
so that for $\lambda\notin [2,\infty)$  the formulas \eqref{R-0-def} and \eqref{R-0-formula} do indeed determine a bounded operator on $L^{2}(\RR)$.

It is instructive to rewrite \eqref{R-0-formula} in terms of the solutions of the equation
\begin{equation} \label{H-0-eqn}
\psi(x+ib-i0,k)+\psi(x-ib+i0,k)=2\cosh(2\pi bk)\psi(x,k)
\end{equation}
for the continuous spectrum eigenvalues of the operator $H_{0}$, that is, in terms of the solutions $f_{\pm}(x,k)=e^{\pm2\pi ikx}$, which are analogues of the Jost solutions in the theory of the one-dimensional Schr\"{o}dinger operator (see \S \ref{Res-Schr}). Namely,
\begin{equation} \label{R-0-kernel-2}
R^{0}_{\lambda}(x-y)=\frac{i}{bC(f_{-},f_{+})(k)}\left(\frac{f_{-}(x,k)f_{+}(y,k)}{1-e^{2\pi (x-y)/b}}+\frac{f_{-}(y,k)f_{+}(x,k)}{1-e^{-2\pi (x-y)/b}} \right),
\end{equation}
where
$$C(f,g)(x,k)=f(x+ib,k)g(x,k)-f(x,k)g(x+ib,k)$$
is the so-called Casorati determinant, which is an analogue of the Wronskian for solutions of the functional-difference equation  
\eqref{H-0-eqn}.  It is periodic function of $x$ with period $ib$, and in the case of the Jost solutions
$C(f_{-},f_{+})(x,k)=2\sinh(2\pi bk)$. 

There is a remarkable similarity between \eqref{R-0-kernel-2} and the formulas \eqref{formula} and \eqref{Res-Scr},  
where instead of the Heaviside function $\theta(x)$ a smoothed analogue of it is involved, namely, the function $\theta_{b}(x)$ defined by the formula\footnote{As noted by A.M. Polyakov, the function $\theta_{b}(x)$, after identification of  $x$ with the energy $\epsilon$ and identification of $\frac{2\pi}{b}$ with the inverse temperature $\frac{1}{kT}$, coincides with the one-particle partition function $\displaystyle{\mathcal{Z}=\left(1-e^{-\frac{\epsilon}{kT}}\right)^{-1}}$ in the Bose-Einstein statistics.}
$$\theta_{b}(x)=\frac{1}{1-e^{-2\pi x/b}}.$$
In this case, the analogue of the relation $\theta'(x)=\delta(x)$ is the formula
$$\frac{1}{ib}\left(\theta_{b}(x-i0)-\theta_{b}(x+i0)\right)=\delta(x)$$
for real $x$, which is equivalent to the Sokhotski-Plemelj formula. The following simple formula also holds
\begin{equation}\label{vp}
\frac{1}{ib}\theta_{b}(x\pm i0)=\frac{1}{ib}\mathrm{v.p.}\,\theta_{b}(x)\mp\frac{1}{2}\delta(x),
\end{equation}
where the distribution $\theta_{b}(x)$ is understood as the Cauchy principal value. From this we obtain for $R^{0}_{\lambda}(x-y)$  the equation
 \begin{equation} \label{delta}
R^{0}_{\lambda}(x-y+ib-i0)+R^{0}_{\lambda}(x-y-ib+i0)-\lambda R^{0}_{\lambda}(x-y)=\delta(x-y).
\end{equation}
Indeed, setting $y=0$ and using  \eqref{H-0-eqn}, \eqref{vp} and the regularity of $R^{0}_{\lambda}(x)$ at $x=0$, we have
\begin{gather*}
R^{0}_{\lambda}(x+ib-i0)+R^{0}_{\lambda}(x-ib+i0)-\lambda R^{0}_{\lambda}(x) \\
=\frac{i}{2b\sinh(2\pi bk)}\left[f_{-}(x+ib,k)\theta_{b}(-x+i0) +f_{+}(x+ib,k)\theta_{b}(x-i0)\right.\\
+\left.f_{-}(x-ib,k)\theta_{b}(-x-i0)+f_{+}(x-ib,k)\theta_{b}(x+i0)\right]\\
-\frac{i}{b}\coth(2\pi bk)(f_{-}(x,k)\theta_{b}(-x)+f_{+}(x,k)\theta_{b}(x))\\
=\frac{1}{4\sinh(2\pi bk)}\left[f_{-}(ib,k)-f_{+}(ib,k)-f_{-}(-ib,k)+f_{+}(-ib,k)\right]\delta(x)=\delta(x).
\end{gather*}
By using the representation \eqref{R-0-kernel-2} and equation \eqref{delta}, it is easy to verify directly that for $\lambda\in\CC\setminus [2,\infty)$ the integral operator \eqref{R-0-def} is the inverse of the operator $H_{0}-\lambda I$ (see \S\S \ref{Sturm-L}--\ref{Schr-one}). 

\subsection{The operator $H$} \label{H-E} Here we consider the equation
\begin{equation}\label{diff-0}
\psi(x+ib-i0)+\psi(x-ib+i0) +e^{2\pi bx}\psi(x)=2\cosh(2\pi bk)\psi(x),
\end{equation}
which is the $q$-analogue of the equation
$$-\psi''+e^{2x}\psi=k^{2}\psi$$
for the Bessel functions. As is well known, the last equation has a solution that is decreasing  as $x\to\infty$, the modified Bessel function of the second kind $K_{ik}(e^{x})$ given by the inverse 
Mellin transform of the product of two gamma functions. The equation \eqref{diff-0} also has a solution that is decreasing  as $x\to\infty$, the Fourier transform of a product involving another wonderful special function, Faddeev's quantum dilogarithm.  
This function was introduced by Faddeev in \cite{F00} and has the integral representation\footnote{The function $\Phi_{b}(z)$ has an interesting history (see \cite{TF15}, where the notation  $\gamma(z)$ was used).} 
\begin{equation} \label{dilog}
\Phi_{b}(z)=\exp\left\{ \frac{1}{4}\int_{-\infty}^{\infty}\frac{e^{2itz}}{\sinh bt\sinh b^{-1}t}\frac{dt}{t}\right\},
\end{equation}
where the contour of integration passes above the singularity at $t=0$. 
The representation \eqref{dilog}  is valid for $|\im z|<c_{b}=\frac{1}{2}(b+b^{-1})$ and defines a meromorphic function with poles $z=-ic_{b} -mib -nib^{-1}$ for integer $m,n\geq 0$, which satisfies the functional equations
\begin{align*}
\Phi_{b}(z+ib) &=(1+q^{-1}e^{-2\pi bz})\Phi_{b}(z),\quad q=e^{\pi ib^{2}},\\
\Phi_{b}(z+ib^{-1}) &=(1+\tilde{q}^{-1}e^{-2\pi b^{-1}z})\Phi_{b}(z),\quad \tilde{q}=e^{\pi i b^{-2}}. 
\end{align*}
Let
$$\hat{\vp}(p,k)=\exp\{-i\beta-\pi ik^{2}-\pi i(p-ic_{b})^{2}\}\Phi_{b}(p-k-ic_{b})\Phi_{b}(p+k-ic_{b}),$$
where $\beta=\dfrac{\pi}{12}(b^{2}+b^{-2})$. 
Using the analytic properties of $\Phi_{b}(z)$ (see, for example, \cite{TF15})), it is easy to verify that the function
$$\vp(x,k)=\int_{-\infty}^{\infty}\hat{\vp}(p,k)e^{2\pi ipx}dp$$
is a solution of \eqref{diff-0},
where the contour of integration passes above the singularities at $p=\pm k$. 
Namely, the following statements hold.
\begin{itemize}
\item[\bf{1.}] For real $k$ the function $\varphi(x,k)$ is an even real-valued function of  $k$, having the asymptotics  
$$\vphi(x,k)=M(k)e^{2\pi ikx}+M(-k)e^{-2\pi ikx} + o(1)$$
as real $x\rightarrow -\infty$,  
where
$$M(k)=\exp\left\{i(\beta+\frac{\pi}{4}) - 2\pi ik(k-ic_{b})\right\}\Phi_{b}(2k-ic_{b}),\quad \overline{M(k)}=M(-k),$$
and
$$\frac{1}{|M(k)|^{2}}=4\sinh (2\pi bk)\sinh(2\pi b^{-1} k).$$
So $\vp(x,k)$ is a scattering solution for \eqref{diff-0}.
\item[\bf{2.}] For real $x$ the function $\varphi(x,k)$ admits analytic continuation into the strip $0<\im k\leq 1/(2b)$ and satisfies the reality condition
$$\overline{\varphi(x,k)}=\varphi(x,-\bar{k}).$$
\item[\bf{3.}] For fixed $k$ in the physical strip, the function $\varphi(x,k)$ extends to an entire function of the complex variable $x$ and satisfies equation \eqref{diff-0}.
\item[\bf{4.}] The following estimates hold:
\begin{equation*}
|\vphi(x,k)|\leq Ce^{-2\pi\im k x},
\end{equation*}
uniformly for $-\infty<x\leq a$,  and
\begin{equation*}
|\vphi(x,k)|\leq C e^{-\pi (b+b^{-1})x},\quad 
|\vphi(x\pm ib,k)|\leq Ce^{\pi(b-b^{-1})x},
\end{equation*}
uniformly for $a\leq x<\infty$.
\end{itemize}

As $x\rightarrow-\infty$, equation \eqref{diff-0} takes on the free form \eqref{H-0-eqn}, so it is natural to assume that \eqref{diff-0} has Jost solutions, that is, solutions $f_{\pm}(x,k)$ with the asymptotics
\begin{equation} \label{f-asym}
f_{\pm}(x,k)=e^{\pm2\pi ikx} + o(1)\quad\text{as}\quad x\rightarrow-\infty.
\end{equation} 
Namely, let
\begin{gather*} 
f_{+}(x,k)=\frac{1}{4 \sinh (2\pi b^{-1} k)M(k)}\\ \times\big(\vphi(x-ib^{-1},k)-\vphi(x+ib^{-1},k)+2\sinh(2\pi b^{-1} k)
\vp(x,k)\big)
\end{gather*}
and $f_{-}(x,k)=f_{+}(x,-k)$. From properties \textbf{1} and \textbf{3} of the function $\vp(x,k)$ we immediately get that for real $x$ the functions $f_{\pm}(x,k)$ are solutions of \eqref{diff-0} and
\begin{equation}\label{scattering}
\vphi(x,k)=M(k)f_{+}(x,k)+M(-k)f_{-}(x,k).
\end{equation}
From the properties of the solution $\vp(x,k)$ listed above it is not difficult to derive the following properties of the Jost solutions.
\begin{itemize} 
\item[$\bf{1^{\prime}.}$] For real $x$ and $k$ the functions $f_{\pm}(x,k)$ have the asymptotics \eqref{f-asym}.
\item[$\bf{2^{\prime}.}$] For real $x$ the functions $f_{\pm}(x,k)$ admit analytic continuation to the physical strip $0< \im k\leq 1/(2b)$ and satisfy the condition
$$\overline{f_{\pm}(x,k)}=f_{\pm}(x,-\bar{k}).$$
\item[$\bf{3^{\prime}.}$] For fixed $k$ in the physical strip, the $f_{\pm}(x,k)$ are entire functions of the variable $x$ and satisfy equation \eqref{diff-0} and condition \eqref{scattering}. Moreover, the asymptotics in  $\bf{1^{\prime}}$ remain valid for $0\leq\im x\leq b$.
\item[$\bf{4^{\prime}.}$] The estimates
$$|f_{\pm}(x,k)|\leq Ce^{\mp2\pi\im k x},$$
hold uniformly for $-\infty<x\leq a$, and
$$|f_{\pm}(x,k)|\leq Ce^{\pi(b^{-1}-b)x},\quad |f_{\pm}(x+ib,k)|\leq Ce^{\pi(b+b^{-1})x},$$
uniformly for  $a\leq x<\infty$.
\end{itemize}

Using these analytic properties and the Phragm\'{e}n-Lindel\"{o}f theorem, one can prove that Casorati determinant of the Jost solutions does not depend on $x$, and therefore
$$C(f_{-}, f_{+})(x,k)= 2\sinh(2\pi bk).$$
Arguing as in case of the free operator $H_{0}$,  from this we get that for $\lambda\in\CC\setminus [2,\infty)$ the integral operator $R_{\lambda}$ acting in $L^{2}(\RR)$ with the symmetric kernel
\begin{gather} \label{R-kernel-2}
R_{\lambda}(x,y)=\frac{i}{2b\sinh(2\pi b k)M(k) }\nonumber\\
\times(f_{-}(x,k)\vp(y,k)\theta_{b}(y-x)+f_{-}(y,k)\vp(x,k)\theta_{b}(x-y)),
\end{gather}
 is the resolvent of the operator $H$. Indeed, since the functions $\vphi(x,k)$ and $f_{-}(x,k)$ satisfy \eqref{diff-0},  we get from \eqref{vp} the equation
\begin{equation} \label{delta2}
R_{\lambda}(x+ib-i0,y)+R_{\lambda}(x-ib+i0,y)+(e^{2\pi bx}-\lambda)R_{\lambda}(x,y)=\delta(x-y),
\end{equation}
and we get from the analytic properties of these functions the estimate
$$|R_{\lambda}(x,y)|\leq Ce^{-2\pi\im k|x-y|},$$
so that for $\lambda\in\CC\setminus [2,\infty)$ the operator $R_{\lambda}$ is bounded on $L^{2}(\RR)$. Using \eqref{delta2} and the identity
$$C(f_{-},\vp)(x,k)=2\sinh(2\pi bk)M(k),$$
we obtain the desired statement $R_{\lambda}=(H-\lambda I)^{-1}$.

Finally, the eigenfunction expansion theorem for the operator $H$ is obtained from \eqref{Stone}. Namely, computing the jump of the resolvent kernel $R_{\lambda}(x,y)$ on the branch cut $[2,\infty)$ using
\eqref{scattering}, we get that the operator $\cU$ given by the formula
$$(\cU\psi)(k)=\int_{-\infty}^{\infty}\psi(x)\varphi(x,k)dx,\quad\psi(x)\in L^{2}(\RR),$$
maps $L^{2}(\RR)$ isometrically onto the Hilbert space $\cH_{0}=L^{2}([0,\infty), |M(k)|^{-2}dk)$, that is,
$$\cU^{*}\cU=I\quad\text{and}\quad \cU\cU^{*}=I_{0},$$
where $I_{0}$ is the identity operator on $\cH_{0}$. Moreover, the operator $\cU H\cU^{-1}$ is the multiplication by the function $2\cosh(2\pi bk)$ operator on $\cH_{0}$, so $H$ has a simple absolutely continuous spectrum filling $[2,\infty)$. As was noted in \cite{TF15},  the eigenfunction expansion theorem for $H$ is a $q$-analogue of the classical Kontorovich-Lebedev transform in  the theory of special functions.

\subsection{Operators for mirror curves} \label{M-C} In \cite{AV06} a remarkable connection was found between the functional-difference operators constructed from Weyl operators $U$ and $V$ and the quantization of algebraic curves that are the images of toric Calabi-Yau threefolds under the mirror symmetry. A typical example of such varieties is the total space of the canonical bundle of a toric del Pezzo surface $S$. The spectral properties of such operators were studied in \cite{GHM16}. In the simplest case, when $S$ is the Hirzebruch surface $S=\mathbb{P}^{1}\times \mathbb{P}^{1}$, we obtain the operator
$$H(\zeta)=U+U^{-1}+V+\zeta V^{-1},$$
where the parameter  $\zeta>0$ plays the role of a physical mass. Thus the operator $H$ considered above corresponds to the massless case $H=H(0)$. When $S=\mathbb{P}(1,m,n)$ is a weighted projective space with $m,n\in\NN$, the corresponding operator has the form
$$H_{m,n} = U + V + q^{-mn}U^{-m}V^{-n}.$$ 

In \cite{LST16} the spectral properties of the self-adjoint operators $H(\zeta)$ and  $H_{m,n}$ in $L^{2}(\RR)$ were investigated. In particular, there it was proved that these operators have a purely discrete spectrum, and an asymptotic expression for the eigenvalues was obtained which implies that $H(\zeta)^{-1}$ and $H_{m,n}^{-1}$ are trace class operators. Furthermore, an analogue of Weyl's asymptotic law was obtained for the eigenvalue counting function
 $N(\lambda)$: it was proved that
$$\lim_{\lambda\to\infty}\frac{N(\lambda)}{\log^{2}\lambda} =\frac{1}{(\pi b)^{2}}$$
for the operator $H(\zeta)$, and
 $$\lim_{\lambda\to\infty}\frac{N(\lambda)}{\log^{2}\lambda} =\frac{c_{m,n}}{(2\pi b)^{2}},\quad\text{where}\quad c_{m,n}=\frac{(m+n+1)^{2}}{2mn}$$
for the operator $H_{m,n}$. Hence for the eigenvalues $\lambda_{k}$ we get that 
 $$\lambda_{k}=e^{\alpha\sqrt{k}}\,(1+o(1)),$$
as  $k\to\infty$, where $\alpha=\pi b$ for $H(\zeta)$ and $\alpha=2\pi b/\sqrt{c_{m,n}}$ for $H_{m,n}$ (see \eqref{l-as} in \S \ref{Res-S-L}). 

A detailed proof of these formulas was given in \cite{LST16}. We note here that it would be quite interesting to obtain more accurate asymptotic formulas for the eigenvalues of the operators $H(\zeta)$ and $H_{m,n}$. As we saw in \S \ref{det-S-L},  it is rather instructive to compare asymptotics of the eigenvalues with the asymptotics of the Fredholm determinants of the operators $H(\zeta)^{-1}$ and $H_{m,n}^{-1}$. In \cite{GHM16}  a remarkable connection was pointed out between these determinants and the enumerative invariants of the corresponding Calabi-Yau manifolds. We leave it to the reader
to reflect on these intriguing connections and associations.

\section{Laplace operator on the fundamental domain of a discrete group on the Lobachevsky plane}\label{Laplace}
Let $\HH=\{z\in\CC : \im z>0\}$ be the Poincar\'{e} model of the Lobachevsky plane with the metric $ds^{2}=\dfrac{|dz|^{2}}{y^{2}}$
and the area form $d\mu(z)=\dfrac{dx\wedge dy}{y^{2}}$, $z=x+iy$.  The group of motions of the Lobachevsky plane is the Lie group $G=\mathrm{PSL}(2,\RR)$, which acts on $\HH$ by fractional-linear transformations,
$$\HH\ni z\mapsto gz=\frac{az+b}{cz+d}\in\HH,\quad\text{where}\quad g=\begin{pmatrix} a & b\\c & d\end{pmatrix}\in G.$$

Denote by $\mathcal{A}$ the Laplace operator of the Poincar\'{e} metric,
\begin{equation}\label{Delta}
\mathcal{A}=-y^{2}\Delta\quad\text{where}\quad \Delta=\frac{\del^{2}}{\del x^{2}} + \frac{\del^{2}}{\del y^{2}},
\end{equation}
defined on the space $C_{0}^{\infty}(\HH)$ of smooth functions with compact support. The operator $\mathcal{A}$ commutes with the action of $G$ on $\HH$ and is essentially self-adjoint on the Hilbert space $\cH_{0}=L^{2}(\HH,d\mu)$. Its closure, a self-adjoint operator $A_{0}=\bar{\mathcal{A}}$, has an absolutely continuous spectrum of infinite multiplicity, filling the interval 
$ [1/4,\infty)$. It is convenient to use a parametrization $\lambda=s(1-s)$, in which the resolvent set  $\CC\setminus [1/4,\infty)$ corresponds to the half-plane $\mathrm{Re}\,s>1/2$. 

The operator $A_{0}$ is invariant under the action of $G$ on $\HH$,  
so its resolvent $R^{0}_{\lambda}=(A_{0} - s(1-s)I)^{-1}$ is the integral operator with kernel
$$R^{0}_{\lambda}(z,z')=r_{0}(z,z';s)=\vp(u(z,z'),s),$$
where
$$u(z,z')=\frac{|z-z'|^{2}}{4yy'}$$ 
is a point-pair invariant in the Lobachevsky geometry ($u(gz,gz')=u(z,z')$ for all $g\in G$), and $\vp(u,s)$ is given by the classical integral
$$\vp(u,s)=\frac{1}{4\pi}\int_{0}^{1}[t(1-t)]^{s-1}(t+u)^{-s}dt $$
and can be expressed explicitly in terms of the hypergeometric function. For fixed $s$ the function $\vp(u,s)$ has asymptotics
\begin{equation}\label{vp-as}
\vp(u,s) =-\frac{1}{4\pi}\log u +O(1) \quad \text{as}\quad u\to 0
\end{equation}
and
\begin{equation}  \label{vp-as-1}
 \vp(u,s)=O(u^{-\sigma})\quad \text{as}\quad  u\to \infty,\quad\text{where}\quad \sigma=\re s.
\end{equation}
Moreover,
\begin{equation} \label{real}
\overline{r_{0}(z,z';s)}=r_{0}(z,z';\bar{s}).
\end{equation}
Using these formulas, it is easy to check directly that if $f\in\cH_{0}$, then
$$(R^{0}_{\lambda}f)(z)=\iint_{\HH}r_{0}(z,z';s)f(z')d\mu(z')\in D(A_{0}),\quad\text{where}\quad\lambda=s(1-s),$$
and  $(A_{0}-\lambda I)R^{0}_{\lambda}f=f$ (see \cite{Fad67}, and also the monographs \cite{Hejhal, Lang}).
\subsection{The resolvent and the eigenfunction expansion}\label{Res-hyp}
Let  $\Gamma$ denote a  Fuchsian group of the first kind, that is, a discrete subgroup of the group $G=\mathrm{PSL}(2,\RR)$ such that the quotient $\Gamma\backslash\HH$ has finite area,
$$\mu(F)=\iint_{F}\frac{dxdy}{y^{2}}<\infty,$$
where $F$ is a fundamental domain of $\Gamma$ in $\HH$. Recall that $F$ is an open subset in $\HH$ such that $\gamma_{1}F\cap\gamma_{2}\bar{F}=\emptyset$ when $\gamma_{1}\neq\gamma_{2}$, and the union $\cup\gamma\bar{F}$ over all $\gamma\in\Gamma$ is $\HH$. Equivalently, a Fuchsian group of the first kind is a discrete subgroup $\Gamma$ of the group $G$ which is finitely generated by hyperbolic  
generators $\alpha_{1},\beta_{1},\dots,\alpha_{g},\beta_{g}$, elliptic generators $\sigma_{1},\dots,\sigma_{m}$ of orders $k_{1},\dots,k_{m}\geq 2$, and parabolic generators $\tau_{1},\dots,\tau_{n}$. They satisfy the relations
$$[\alpha_{1},\beta_{1}]\cdots [\alpha_{g},\beta_{g}]\sigma_{1}\cdots\sigma_{m}\tau_{1}\cdots\tau_{n}=1\quad \text{and}\quad \sigma_{1}^{k_{1}}=\cdots=\sigma_{m}^{k_{m}}=1,$$
where $[\alpha,\beta]=\alpha\beta\alpha^{-1}\beta^{-1}$, and the condition
$$\chi(\Gamma)=2-2g - n -\sum_{j=1}^{m}\left(1-\frac{1}{k_{j}}\right)<0,$$
wherein $\mu(F)=-2\pi\chi(\Gamma)$. In case $n>0$ the closure $\bar{F}$ of the fundamental domain $F$ is non-compact in $\HH$ and contains $n$ cusps, fixed points of the parabolic transformations  $\tau_{1},\dots,\tau_{n}$,  lying on $\RR\cup\{\infty\}$. The simplest example is known to Gauss (see \cite{K}) fundamental domain $F$ of the modular group $\mathrm{PSL}(2,\ZZ)$,
$$F=\left\{z\in \HH : |x|<\tfrac{1}{2}\;\;\text{и}\;\; x^{2}+y^{2}>1\right\}.$$ 
We have $ \mathrm{PSL}(2,\ZZ)\bk\HH\simeq \tilde{F}$, where $ \tilde{F}$ is the so-called modular figure,
$$\tilde{F}=\left\{z\in\HH: -\tfrac{1}{2}< x< 0 \;\;\text{и}\;\; x^{2}+y^{2}>1\right\}\cup\left\{z\in\HH: 0\leq x\leq \tfrac{1}{2} \;\;\text{и}\;\; x^{2}+y^{2}
\geq1\right\}.$$

A measurable function $f$ on $\HH$ is called $\Gamma$-automorphic, if
$$f(\gamma z)=f(z)$$
for all  $z\in\HH$ and $\gamma\in\Gamma$.  Let  $\cH=L^{2}(F,d\mu)$ be the Hilbert space of $\Gamma$-automorphic functions that are 
square-integrable on  $F$ with respect to the measure $d\mu$, with
$$\Vert f\Vert^{2}=\iint_{F}|f(z)|^{2}\frac{dxdy}{y^{2}}<\infty.$$
It is not difficult to show that the differential expression \eqref{Delta}, defined on the space $C_{0}^{\infty}(F)$  of smooth functions on  $F$ with compact support, is an essentially self-adjoint operator in $\cH$. Denote its closure by $A$.

Since the Laplace operator in the spaces $\cH_{0}$ and $\cH$ is given by the same differential expression \eqref{Delta}, it is natural to assume that the resolvent of $A$ in $\cH$, that is, the operator
$$R_{\lambda}=(A - s(1-s)I)^{-1},\quad\text{where}\quad\lambda=s(1-s),$$ 
is still the integral operator with the integral kernel $R_{\lambda}(z,z')=r(z,z';s)$ obtained from the kernel $r_{0}(z,z';s)$ by the classical method of images. Using a simple criterion for the convergence over a discrete group  (see \cite{Fad67,Lang}) and the estimate \eqref{vp-as-1}, we can easily prove that if $z\neq \gamma z'$ for all $\gamma\in\Gamma$, then for $\sigma>1$ the series 
\begin{equation} \label{series}
r(z,z';s)=\sum_{\gamma\in\Gamma}r_{0}(z,\gamma z';s),
\end{equation}
is absolutely convergent, uniformly with respect to $z, z'$ on every compact subset, and it satisfies reality condition \eqref{real}. The further analysis depends essentially  on whether the closure $\bar{F}$ of the fundamental domain in  $\HH$ is compact (the case $n=0$) or non-compact (the case $n\geq 1$). 

 The case $n=0$ is elementary. Indeed, for $\sigma>1$ the kernel $r(z,z';s)$  has a weak singularity on the diagonal in $F\times F$ and defines a compact operator on $\cH$, the resolvent $R_{\lambda}$ of the operator $A$. The eigenfunction expansion theorem immediately follows from the first Hilbert identity
  \eqref{Hilbert-1},
$$R_{\lambda}-R_{\mu}=(\lambda-\mu)R_{\mu}R_{\lambda},$$
and the Hilbert-Schmidt decomposition for the compact operator.

Indeed, choose $\vk>1$ and put $R=R_{\mu}$, where $\mu=\vk(1-\vk)<0$. The self-adjoint compact operator $R$ is positive, so
$$R=\sum_{n=1}^{\infty}\mu_{n}\PPP_{n},$$
where $\PPP_{n}$ are orthogonal projection operators in $\cH$ onto finite-dimensional invariant subspaces of the operator $R$ corresponding to the eigenvalues $\mu_{n}>0$, and $\mu_{1}=-\mu^{-1}$.
Here
$$\sum_{n=1}^{\infty}\PPP_{n}=I\quad\text{and}\quad\lim_{n\to\infty}\mu_{n}=0.$$
Rewriting the Hilbert identity as equation for $R_{\lambda}$,
\begin{equation}\label{H-2}
(I-(\lambda-\mu)R)R_{\lambda}=R,
\end{equation}
we get the eigenfunction expansion theorem for the operator $A$:
$$R_{\lambda}=\sum_{n=1}^{\infty}\frac{\PPP_{n}}{\lambda_{n}-\lambda},\quad\text{where}\quad\lambda_{n}=\mu+\frac{1}{\mu_{n}}.$$

In case  $n\geq 1$ the derivation of the eigenfunction expansion theorem for $A$ is rather complicated. Namely, the spectrum of the Laplace operator now consists of the $n$-fold absolutely continuous spectrum $[1/4,\infty)$ and finite multiplicity eigenvalues lying on $0\leq\lambda<\infty$ without accumulation points on a finite interval. Moreover, the so-called Eisenstein-Maass series, defined by the series over the cosets of $\Gamma$ which are absolutely convergent for $\sigma=\re s>1$ admit meromorphic continuation to the whole complex $s$-plane, with poles for $\sigma<1/2$ and on the interval $[1/2,1]$, and the eigenfunctions of the continuous spectrum of $A$ are given by analytic continuation of the Eisenstein-Maass series to the line $\sigma=1/2$. These results were announced\footnote{A proof based on  potential theory was presented in Selberg's then unpublished 1954 lectures at the University of G\"{o}ttingen.} in Selberg's famous paper \cite{Selberg}, and were first proved by Faddeev \cite{Fad67}.  The monograph \cite{Lang} of Lang is devoted to a detailed presentation of Faddeev's method. 

It is easy to show \cite{Fad67} that the series \eqref{series}  still converges for $\sigma>1$, and for $\sigma>2$ it is the integral kernel of the resolvent  $R_{\lambda}$ of $A$. 
Nevertheless, equation \eqref{H-2} is no longer suitable for an investigation of the resolvent $R_{\lambda}$ for all $\lambda=s(1-s)\in \CC\setminus [1/4,\infty)$. The fact is that for $\mu=\vk(1-\vk)$ with $\vk>2$ the operator $R=R_{\mu}$ is no longer compact, but rather has an absolutely continuous spectrum. A remarkable observation made by Faddeev in \cite{Fad67} is that the main part of $R$ generating this spectrum 
can be identified and explicitly inverted!  The paper \cite{Fad67} is based on the virtuoso use of the resolvent technique, the spectral theory of Sturm-Liouville operators,  and the Fredholm theory. Here we present only the main steps of the algebraic scheme of calculations; a detailed derivation of all necessary estimates can be found in \cite{Fad67,Lang}.

In particular, for simplicity consider the case of one cusp\footnote{The case of several cusps is considered in the same way.}
$ i \infty $ and choose a fundamental domain $ F $ in the form
$$F=F_{0}\cup F_{1},$$
where $\bar{F}_{0}$ is compact and $\bar{F_{1}}$ is a strip $\{z=x+iy : -\frac{1}{2}\leq x\leq \frac{1}{2}, \;y\geq a\}$ for some $a>0$. Denote by $P_{0}$ and $P_{1}=I-P_{0}$ orthogonal projection operators on  $\cH$ corresponding to multiplication by the indicator functions of the regions $F_{0}$ and $F_{1}$, and write the operator  $R$ as
$$R=R_{00}+R_{10}+R_{01}+R_{11},$$
where $R_{00}=P_{0}RP_{0}$, $R_{01}=P_{0}RP_{1}$, $R_{10}=P_{1}RP_{0}$ and $R_{11}=P_{1}RP_{1}$. 

By using \eqref{series} it is not difficult to show that for $\varkappa>2$ the operators  $R_{00}, R_{01}$ and $R_{10}$ are compact. It follows from the representation \eqref{series}, that the `cusp part' of $R$, that is, the operator $R_{11}=P_{1}RP_{1}$, is the integral operator with kernel
$$R_{11}(z,z')=\sum_{\gamma\in\Gamma_{\infty}}R_{\vk}^{0}(z,\gamma z'),$$
where $\Gamma_{\infty}=\left\{\begin{pmatrix} 1 & n\\0 & 1\end{pmatrix}, n\in\ZZ\right\}$ is the stabilizer of the cusp $i\infty$ in the group $\Gamma$. Rewriting this formula as
$$R_{11}(z,z')=\sum_{n=-\infty}^{\infty}\vp\left(\frac{|z-z'-n|^{2}}{4yy'},\vk\right),$$
we see that the kernel $R_{11}(z,z')$ is an even periodic function of the variable $x-x'$ with period $1$, and it can be expanded in a Fourier series
\begin{equation}\label{F}
R_{11}(z,z')=t_{0}(y,y')+\sum_{m=1}^{\infty}t_{k}(y,y')\cos (2\pi m(x-x')).
\end{equation}
The constant term, that is, the function $t_{0}(y,y')$, is easily computable:
$$t_{0}(y,y')=\frac{1}{2\vk-1}\begin{cases}y^{\vk}y'^{1-\vk}, & y\leq y',\\
y^{1-\vk}y'^{\vk}, & y> y', \end{cases}$$
and for the functions $t_{k}(y,y')$ it is not difficult to obtain expressions in terms of the modified Bessel functions. 

The operator $R_{11}$ acts in the Hilbert space $P_{1}\cH$, which has a natural subspace $L^{2}([a, \infty); y^{-2}dy)$ consisting of functions independent of $x$,  
and the orthogonal projection operator $P$ from $P_{1}\cH$ onto $L^{2}([a, \infty); y^{-2}dy)$ is given by the integration:
\begin{equation}\label{c-term}
f(z)\mapsto P(f)(y)=\int_{-\frac{1}{2}}^{\frac{1}{2}}f(x+iy)dx,\quad y\geq a.
\end{equation}
We write the kernel $R_{11}(z,z')$ in the form
$$R_{11}=T+ R_{11}',$$
where $T=PR_{11}P$ is the integral operator in $L^{2}([a, \infty); y^{-2}dy)$ with kernel $t_{0}(y,y')$, 
and $R'_{11}$ is the integral operator with kernel $R_{11}(z,z')-t_{0}(y,y')$. Using the Fourier expansion \eqref{F} and standard estimates for the modified Bessel functions, we can easily prove (for details, see \cite{Lang}), that the operator $R'_{11}$ is compact.
By recalling formulas in the \S \ref{S-L}, it is easy to verify that
$$T=(B-\vk(1-\vk)I)^{-1},$$ 
where $B$ is a self-adjoint operator acting in  $L^{2}([a, \infty); y^{-2}dy)$ and given by the differential expression
$\mathcal{B}\vp=-y^{2} d^{2}\vp/dy^{2}$ and the boundary condition\footnote{Here we correct the typing error in \cite{Fad67} after (3.7) and also in the corresponding place in \cite[\S 3.1]{T17}.}
$$\vk\varphi(a)=a\varphi'(a).$$
For $\sigma>1/2$ the resolvent $Q_{\lambda}=(B-s(1-s)I)^{-1}$ is the integral operator with kernel
$$q(y,y';s)=\frac{1}{2s-1}\left(\vp(y,s)y'^{1-s}\theta(y'-y)+y^{1-s}\vp(y',s)\theta(y-y')\right),$$
where
$$\vp(y,s)=y^{s}+a^{2s-1}\frac{s-\vk}{s+\vk-1}y^{1-s}.$$
For $\sigma=1/2$ the functions $\vp(y,s)$ form a complete system of the continuous spectrum eigenfunctions for the operator $B$ in the space $L^{2}([a, \infty); y^{-2}dy)$.

Summing up, for $\vk>2$ we have
$$R=T+V,$$ 
where $V$ is a compact operator. Therefore, $A$ can be regarded as a perturbation of $B$ with the same absolutely continuous spectrum! Namely, we now write 
 \eqref{H-2} in the form 
\begin{equation}\label{H-3}
(I-(\lambda-\mu)T)R_{\lambda}= R + (\lambda-\mu)VR_{\lambda},
\end{equation}
where $\lambda=s(1-s)$ and  $\mu=\vk(1-\vk)$. It follows from the first Hilbert identity for $B$ that
$$I-(\lambda-\mu)T=\left(I+(\lambda-\mu)Q_{\lambda}\right)^{-1},$$
so \eqref{H-3} can be rewritten as
$$R_{\lambda}=(I+(\lambda-\mu)Q_{\lambda})(T+V)+(\lambda-\mu)(I+(\lambda-\mu)Q_{\lambda})VR_{\lambda},$$
or
\begin{equation}\label{H-4}
R_{\lambda}=Q_{\lambda}+(I+(\lambda-\mu)Q_{\lambda})V+(\lambda-\mu)(I+(\lambda-\mu)Q_{\lambda})VR_{\lambda}
\end{equation}
if one uses the Hilbert identity once again. Putting
\begin{equation}\label{R-U}
R_{\lambda}=Q_{\lambda}+\left(I+(\lambda-\mu)Q_{\lambda}\right)U_{\lambda}\left(I+(\lambda-\mu)Q_{\lambda}\right),
\end{equation}
we obtain for $U_{\lambda}$ the equation
\begin{equation}\label{H-5}
U_{\lambda}=V+H_{\lambda}U_{\lambda},\quad\text{where}\quad H_{\lambda}=(\lambda-\mu)V\!\left(I+(\lambda-\mu)Q_{\lambda}\right).
\end{equation}

Equation \eqref{H-5}, Faddeev's equation in the theory of automorphic functions, has the following remarkable properties \cite{Fad67}. 
\begin{itemize}
\item[1)]  The operator $H_{\lambda}$ is a Fredholm operator in the Banach space
$\frak{B}$ of continuous functions $f(z)$  on $F$ with the norm
$$\Vert f\Vert_{\frak{B}} =\sup_{z\in F_{0}}|f(z)|+\sup_{z\in F_{1}}y|f(z)|$$
and depends analytically on $s$ in the strip $0<\sigma< 2$. 
\item[2)] The singular points of the operator $I-H_{\lambda}$, that is,  the values of $s$ for which the homogeneous equation
$$v=H_{\lambda}v$$
has a nontrivial solution in the space $\frak{B}$, are discrete in the strip $0<\sigma< 2$.
\item[3)]  The singular points with $\sigma\geq 1/2$, $s\neq  1/2$, correspond to non-negative eigenvalues $\lambda=s(1-s)$ of $A$ of finite multiplicity, so that $\sigma=1/2$ or $1/2<s\leq 1$.
The corresponding eigenfunctions $\psi\in\cH$ have the form
$$\psi=(I+(\lambda-\mu)Q_{\lambda})v,$$
where $v\in\frak{B}$ is a solution of the homogeneous equation. Eigenfunctions corresponding to the case $\sigma=1/2$ are cusp forms\footnote{In general, the space $\cH^{(0)}$ of cusp forms  is an invariant subspace of $\cH$ consisting of functions with zero integrals over all horocycles in $\Gamma\bk\HH$. It is not difficult to show \cite{GPS} that the spectrum of $A$ in $\cH^{(0)}$ is discrete.}, that is, $P(\psi)(y)=0$ for all  $y>0$.
\item[4)] The resolvent kernel $r(z,z';s)$ of $A$ for fixed $z\neq z'$ and $\sigma>1$ admits meromorphic continuation to the strip $0<\sigma <2$, with discrete poles of finite multiplicity. For $\sigma\geq 1/2$ these poles lie only on the line $\sigma= 1/2$ and on the interval $ 1/2\leq s\leq 1$ and are simple,  with the possible exception for $s= 1/2$. 
\item[5)] The resolvent $(A-\lambda I)^{-1}$ of $A$, where $\lambda=s(1-s)\in \CC\setminus [1/4,\infty)$ with non-singular $s$ and $\sigma> 1/2$, is the operator $R_{\lambda}$ in \eqref{R-U}, constructed from the solution $U_{\lambda}$ of the Faddeev's equation \eqref{H-5}. The operator $R_{\lambda}$ is an integral operator in $\cH$ with the integral kernel $r(z,z';s)$.
\end{itemize}

Equation \eqref{H-5} is also used for analytic continuation of the continuous spectrum eigenfunctions of $A$. In particular, consider the decomposition $F=F_{0}\cup F_{1}$ and define the function $\psi(z,s)$ on $F$ by $\psi(z,s)=\vp(z,s)$ for $z\in F_{1}$ and $\psi(z,s)=0$ for $z\in F_{0}$.  Clearly, if  $a$ is large enough, then $\psi(z,s)$ determines a piece-wise smooth  $\Gamma$-automorphic function on $\HH$. We put
$$\Psi(z,s)=(I+(\lambda-\mu)(I+(\lambda-\mu)Q_{\lambda})U_{\lambda})\psi(z,s)$$
and list the properties of $\Psi(z,s)$ \cite{Fad67}.
\begin{itemize}
\item[(i)] For fixed $z$ the function $\Psi(z,s)$ is analytic in the strip $0<\sigma<2$, 
except for the singular points for which $\sigma< 1/2$ or
$ 1/2\leq s\leq 1$, and $\Psi(z,s)$ is analytic In a neighborhood of the line $\sigma= 1/2$, with the possible exception of the point $s= 1/2$.
\item[(ii)] For non-singular $s$ in the strip $0<\sigma<2$  the function $\Psi(z,s)$ is a smooth $\Gamma$-automorphic function on $\HH$, satisfying the equation
\begin{equation}\label{E-M-eq}
-y^{2}\left(\frac{\del^{2}}{\del x^{2}}+\frac{\del^{2}}{\del y^{2}}\right)\Psi(z,s)=s(1-s)\Psi(z,s).
\end{equation}
\end{itemize}

For  $\sigma>1$ the solution of  \eqref{E-M-eq} can be found `explicitly' as the Eisenstein-Maass series $E(z,s)$:
\begin{equation}\label{E-M}
E(z,s)=\sum_{\gamma\in\Gamma_{\infty}\backslash\Gamma}y^{s}(\gamma z).
\end{equation}
Namely, it is easy to show that for $\sigma>1$ the series converges absolutely and uniformly on compact subsets of $\HH$, and defines a $\Gamma$-automorphic function satisfying equation \eqref{E-M-eq}.
For $\sigma>1$ it is not difficult to prove the equality $\Psi(z,s)=E(z,s)$, which gives a meromorphic continuation of $E(z,s)$ to the strip $0<\sigma\leq 1$,  and on the line $\sigma= 1/2$ the function $E(z,s)$ has no singularities except, possibly, the point $s= 1/2$.

Finally, the eigenfunction expansion theorem for $A$ is obtained from the above results using  \eqref{Stone}. The reader can find detailed proofs in Faddeev's paper \cite{Fad67}, the indicated book by Lang \cite{Lang} and  Venkov's   monograph \cite{Venkov}, which generalizes Faddeev's method to vector-valued functions.  The characteristic determinant of the operator $A$ 
is defined using an appropriate regularization of the formula  \eqref{Res-det} and M.G. Krein method of the spectral shift function.  
Moreover, the characteristic determinant of $A$ is expressed in terms of the Selberg zeta function of the Fuchsian group  $\Gamma$, and the calculation  of the regularized trace in \eqref{Res-det} reduces to the famous Selberg trace formula! We refer the reader to \cite{VKF73} for details of these nontrivial calculations. This concludes our exposition of Faddeev's method.

As an interesting example, consider the case of the modular group $\Gamma=\mathrm{PSL}(2,\ZZ)$. The corresponding  Eisenstein-Maass series $E(z,s)$ admits a simple expression in terms of the Epstein zeta function of the positive-definite binary quadratic form $Q(m,n)=am^{2}+bmn+cn^{2}$ of the discriminant $b^{2}-4ac=d<0$, where $a=1$, $b=-2x$ and $c=x^{2}+y^{2}$, so $d=-4y^{2}$. Furthermore, $z=x+iy\in\HH$ is the root of the quadratic form  $Q$,
$$z=\dfrac{-b+\sqrt{d}}{2a}.$$
In particular, from  \eqref{E-M} we easily obtain 
\begin{equation} \label{Ep}
2\zeta(2s)y^{-s}E(z,s)=\sideset{}{'}\sum_{m,n=-\infty}^{\infty}\frac{1}{Q(m,n)^{s}},
\end{equation}
where $\zeta(s)$ is the Riemann zeta function,  and the prime on the summation sign indicates that the term with $m=n=0$ is omitted. 
The Fourier series expansion of the function $E(z,s)$ is given by the beautiful formula\footnote{\eqref{E-mod}  is sometimes called the Selberg-Chowla formula.}
\begin{equation}\label{E-mod}
E(z,s)=y^{s}+c(s)y^{1-s}+\frac{4\sqrt{y}}{\xi(2s)}\sum_{n=1}^{\infty}\sigma_{1-2s}(n)n^{s-\frac{1}{2}}K_{s-\frac{1}{2}}(2\pi ny)\cos( 2\pi nx),
\end{equation}
where $K_{s}(y)$ is the modified Bessel function of the second kind and
$$\sigma_{s}(n)=\sum_{d|n}d^{s},\quad c(s)=\frac{\xi(2s-1)}{\xi(2s)},\quad \xi(s)=\pi^{-\frac{s}{2}}\Gamma\left(\frac{s}{2}\right)\zeta(s).$$

It follows from Faddeev's  method that for fixed $z$ the Eisenstein-Maass series $E(z,s)$ is a holomorphic function on the `physical sheet' $\sigma=\re s> 1/2$ and is regular on the line $\sigma= 1/2$. From here it immediately follows that the zeta function $\zeta(s)$ does not vanish on the line $\sigma=1$, which implies the asymptotic law of primes! However, this method does not 
give any information about the poles of $E(z,s)$ on the `non-physical sheet' $\sigma < 1/2$. One can only say that the non-tirival zeros of $\zeta(s)$ are related to the so-called resonances of the Laplace operator on the modular figure.

\subsection{Pseudo-cusp forms and zeros of $L$-series} \label{Pseudo} In 1977, in  H. Haas's diploma work at the University of Heidelberg under the direction of H. Neuenh\"{o}ffer, several of the first eigenvalues of the discrete spectrum
of the Laplace operator on the modular figure were calculated.  Stark and Hejhal soon noticed that if one writes
 $\lambda_{k}=\frac{1}{4}+t_{k}^{2}$, then the values $s_{k}=\frac{1}{2}+it_{k}$ correspond to the first non-trivial zeros of the Riemann zeta function and the Dirichlet $L$-series $L(s,\chi)$ with the quadratic character modulo 3! This unexpected observation caused a sensation and was actively discussed in correspondence between Cartier and Weil in 1979 \cite{Cartier}, as well as by  Venkov, A.I. Vinogradov, Faddeev and author in the Leningrad branch of V.A. Steklov Mathematical Institute of the USSR Academy of Sciences. Hejhal has decided to  verify Haas's calculations and did not find these zeros among the eigenvalues of the Laplace operator.  

What was the reason for this discrepancy? As Hejhal explained in \cite{Hejhal0}, Haas was using the standard collocation method for the Neumann problem on the modular figure, and he did not notice the appearance of a logarithmic singularity at the corners of the modular figure
$z=\rho$ and $z=i$, where  $\rho=(1+\sqrt{-3})/2$. In particular, the function $f(z)=r(z,z_{0};s)$ for $z\neq z_{0}$ satisfies the equation
\begin{equation} \label{eq}
\mathcal{A}f=\lambda f,\quad\text{where}\quad \lambda=s(1-s),
\end{equation}
 and if $z_{0}=\rho$ or $z_{0}=i$, then with discrete approximation it is easy to miss logarithmic singularity \eqref{vp-as} at $z\to z_{0}$. For $y\to\infty$ and fixed  $z_{0}$,  the resolvent kernel has the asymptotics \cite{Fay}
 $$r(z,z_{0};s)=\frac{y^{1-s}}{2s-1}E(z_{0},s)+O(e^{-2\pi y}),$$
and therefore if $E(z_{0},s)=0$, then $f(z)\in L^{2}(F,d\mu)$. It is remarkable that for $z_{0}=i$ and $z_{0}=\rho$ the function $\zeta(2s)E(z_{0},s)$ is proportional to the Dedekind zeta function of the imaginary quadratic fields $\QQ(\sqrt{-1})$ and $\QQ(\sqrt{-3})$, so $\lambda=s(1-s)$ can be expressed in terms of zeros of $\zeta(s)$ and corresponding $L$-series. However, $f(z)$ is not an eigenfunction of the Laplace operator, since it is not a cusp form. Namely, the condition
 $$\int_{-\frac{1}{2}}^{\frac{1}{2}}f(x+iy)dx=0$$
 holds only for $y>\im z_{0}$; such functions are called pseudo-cusp forms.
Moreover, for $z_{0}$ one can take any point on the modular figure, for example $z_{0}=\sqrt{-5}$.
It is well-known that the function $\zeta(2s)E(\sqrt{-5},s)$ has zeros outside the line $\re s=1/2$, so corresponding values of $\lambda$ will not even be real. 
Thus, the pseudo-cusp forms have no relation to the discrete spectrum of the Laplace operator, and equation \eqref{eq} does not impose any restrictions on $\lambda$.
 
Indeed, if $f\in L^{2}(F,d\mu)$ satisfies \eqref{eq} and $f\in D(A)$,  then from the self-adjointness of the operator $A$ we obtain
$$(\lambda-\bar\lambda)\Vert f\Vert^{2}=(Af,f)-(f,Af)=0.$$
However, although $f(z)=r(z,z_{0};s)\in L^{2}(F,d\mu)$ when $\zeta(2s)E(z_{0},s)=0$, $f\notin D(A)$ and the integral $(Af,f)$ is divergent, hence the previous argument does not apply. Specifically, 
$$(Af,f)=\iint_{F}\mathcal{A}f(z)\overline{f(z)}d\mu(z)=\lambda\Vert f\Vert^{2}+r(z_{0},z_{0};s),$$
where the second term is obviously divergent. By using the reality condition \eqref{real}, the difference $(Af,f)-(f,Af)$ can be defined as a limit, which one easily computes via the first Hilbert identity:
$$\lim_{z\to z_{0}}(r(z,z_{0};s)-r(z,z_{0},\bar{s}))=(\lambda-\bar\lambda)\Vert f\Vert^{2}.$$
The last formula does not give any restriction on $\lambda=s(1-s)$, except for the assumption that $\zeta(2s)E(z_{0},s)=0$.
   
Nevertheless, it makes sense to consider Hilbert spaces of pseudo-cusp forms
$$\cH_{a}=\left\{f\in\cH: \int_{-\frac{1}{2}}^{\frac{1}{2}}f(x+iy)dx=0\quad\text{for}\quad y\geq a\right\},$$
for a fixed $a>0$.
In particular, denote by $\Delta_{a}$ the Friedrichs extension of the operator $\Delta$ restricted to the subspace of smooth functions with compact support in $\cH_{a}$. 
Lax and Philips proved \cite{LP} that the self-adjoint operator $\Delta_{a}$ in $\cH_{a}$ has a purely discrete spectrum, which was studied by Colin de Verdi\`{e}re \cite{CV1,CV2}.
Furthermore, it was suggested in \cite{CV2} that the discrete spectrum of the operator $\Delta_{a}$ for $a=\sqrt{3}/2$ is related to zeros of the Dedekind zeta function of the imaginary quadratic field $\QQ(\sqrt{-3})$.
 
\subsection{Heegner points and Linnik asymptotics} \label{Heegner} The formula \eqref{Ep} provides an explicit expression for the Dedekind zeta function $\zeta_{K}(s)$ of the imaginary quadratic field
$K=\QQ(\sqrt{d})$ of the fundamental discriminant $d<0$ in terms of the Eisenstein-Maass series. As is well known  (e.g. \cite{BS}), the ideal class group of the field $K$ is isomorphic to the group of classes of properly equivalent primitive, positive-definite, integral binary quadratic forms with discriminant $d$. Each such quadratic form can be written as $Q(m,n)=am^{2}+bmn+cn^{2}$ with integer coefficients $a$, $b$, $c$, satisfying
$$a>0,\quad (a,b,c)=1\quad\text{and}\quad b^{2}-4ac=d.$$
The root $z_{Q}$ of the quadratic form $Q$ is given by
\begin{equation*}
z_{Q}=\frac{-b+\sqrt{d}}{2a}\in\HH,
\end{equation*}
and the proper equivalence class of $Q$ is completely determined by the condition $z_{Q}\in\tilde{F}$, where $\tilde{F}$ is the modular figure.  
Points $z_{Q}\in \tilde{F}$ are called Heegner points for the discriminant $d$.  
From this we obtain\footnote{For the details see  \cite{VT-2}, for example, where zeta functions of orders in imaginary quadratic fields are also considered.} 
\begin{equation} \label{Hecke}
\zeta_{K}(s)=\frac{2}{w_{d}}\!\left(\frac{|d|}{4}\right)^{\! -s/2}\!\!\zeta(2s)\sum_{i=1}^{h(d)}E(z_{i},s),
\end{equation}
where $h(d)$ is the ideal class number of the field $K$, $w_{d}$ is the number of units in $K$, and $z_{i}$ run over all Heegner points $z_{Q}$ of the discriminant $d$.

Siegel's celebrated theorem \cite{Siegel} states that for every $0<\vep<1/2$
$$h(d)>c(\vep)|d|^{\frac{1}{2}-\vep}$$
with a non-effective constant $c(\vep)>0$. Before the classical paper \cite{Siegel}, the only known result was Hecke's theorem that the generalized Riemann hypothesis for all $L$-series with quadratic characters implies that $h(d)\to\infty$ as $d\to-\infty$. Surprisingly, in 1933 Deuring \cite{Deuring} proved an unexpected result that the condition $h(d)=1$ for infinitely many negative fundamental discriminants implies the Riemann hypothesis!

Indeed, if $h(d)=1$, then 
$$z_{Q}=\begin{cases} \dfrac{1+\sqrt{d}}{2} & \text{ if $d=4D$},\\
\sqrt{D}  & \text{if $d=4D$ and $D \equiv 2,3\!\!\!\!\pmod 4$}.
\end{cases}$$
In the latter case, we get immediately from \eqref{E-mod} and \eqref{Hecke}  that for such $d$
\begin{equation}\label{one-class}
\zeta_{K}(s)=
|D|^{-s/2}\zeta(2s) E(\sqrt{D},s)=\zeta(2s)(1+ c(s)|D|^{\frac{1}{2}-s})+O(e^{-2\pi |D|}). 
\end{equation}
Suppose now that $\zeta(\rho)=0$ and $\re \rho>1/2$.
Because $\zeta_{K}(s)=\zeta(s)L(s,\chi_{d})$, where $\chi_{d}$ is a quadratic character modulo $d$ given by the Kronecker symbol, by passing to the limit $d\to-\infty$ in \eqref{one-class} we obtain  $\zeta(2\rho)=0$, a contradiction. The case  $d\equiv 1\pmod 4$ is considered similarly.

Mordell soon  \cite{Mordell} generalized Deuring's result and proved that if the class number takes a fixed value for infinitely many fundamental negative discriminants, then the Riemann hypothesis is true.  His proof also uses the formulas \eqref{E-mod} and \eqref{Hecke}.  Finally, Heilbronn \cite{H} used the same assumption to deduce the generalized Riemann hypothesis for all Dirichlet $L$-series with quadratic characters. From this and the aforementioned theorem of Hecke it follows that  $h(d)\to\infty$ as $d\to-\infty$. 

However, in the same year of 1934, Siegel proved his famous theorem, which naturally moved the Deuring--Mordell--Heilbronn method to the background. It was only in the 1960s that some of their arguments were used to solve the celebrated tenth discriminant problem of Gauss, in which Heegner points played a prominent role.

Since $h(d)\to\infty$ as $d\to-\infty$, the question arises of the distribution of the Heegner points on the modular figure. The equivalent problem of the distribution of the integer points on the reduction domain of a two-sheeted hyperboloid $b^{2}-4ac=d<0$ was solved by Linnik \cite{Linnik} using his ergodic method. Here the condition  
$$\left(\frac{d}{p}\right)=1$$
was also assumed for some prime  $p$, 
where $(\frac{n}{p})$ is the Legendre symbol.  In particular,  Linnik  proved that as $d\to-\infty$, the Heegner points are uniformly distributed on the modular figure with respect to the measure 
$$d\mu^{*}=\frac{3}{\pi}\frac{dxdy}{y^{2}},\quad\text{so that}\quad\mu^{*}(F)=1.$$ 
In \cite{VT-1} the uniform distribution was proved on average over $d$, that is, for those values of $d$ for which so-called Sali\'{e} sums admit a good estimate. Finally, Duke \cite{Duke} proved the uniform distribution of Heegner points as $d\to-\infty$ using a non-trivial estimate for the Fourier coefficients of modular forms of half-integer weight, obtained by Iwaniec \cite{Ivan}. Specifically, let $\Omega$ be a convex domain with piece-wise smooth boundary on the modular figure,  and let  $N(\Omega)$ be the number of Heegner points in  $\Omega$. Then the Linnik asymptotic expression
\begin{equation}\label{uniform}
\frac{N(\Omega)}{h(d)}
=\mu^{*}(\Omega) + O(|d|^{-\delta}),
\end{equation}
is valid for some $\delta>0$ (possibly depending on $\Omega$).

We now return to representation \eqref{Hecke} for $\zeta_{K}(s)$ and, as proposed in \cite{VT-1},  we use the uniform distribution of Heegner points on $\tilde{F}_{d}$ --- the modular figure $\tilde{F}$ with the restriction $\im z\leq \sqrt{|d|}/2$. More precisely, assuming that the $\delta>0$ in \eqref{uniform} does not depend on the domain $\Omega$, we replace the sum in \eqref{Hecke} by an integral!  As a result, as $d\to-\infty$ we get that
\begin{equation} \label{zeta-int}
\zeta_{K}(s)=\left(\frac{|d|}{4}\right)^{\! -s/2}\!\!\zeta(2s)h(d)\iint_{\tilde{F}_{d}}E(z,s)d\mu^{*}(z) +O(h(d)|d|^{-\delta-\sigma/2}),
\end{equation}
where $\sigma=\re s$. The integral in \eqref{zeta-int} can be evaluated explicitly. Namely, by using \eqref{E-M-eq}, the integral Green's formula, the invariance of $E(z,s)$ with respect to the modular group, and the Fourier expansion \eqref{E-mod}, we obtain
\begin{align*}
\iint_{\tilde{F}_{d}}E(z,s)d\mu(z) &=\frac{1}{s(s-1)}\iint_{\tilde{F}_{d}}\Delta E(z,s)dxdy\\
& =\frac{1}{s(s-1)}\int_{-\frac{1}{2}}^{\frac{1}{2}}\left.\frac{\del E}{\del y}(z,s)\right|_{y=\tfrac{\sqrt{|d|}}{2}}dx\\
&=\frac{1}{s-1}\left(\frac{|d|}{4}\right)^{\!(s-1)/2}\!\!-\;\;\frac{c(s)}{s}\left(\frac{ |d|}{4}\right)^{\!-s/2}.
\end{align*}
Thus, for fixed $s$ we have
\begin{equation} \label{zeta-formula}
\zeta_{K}(s)=\frac{6h(d)}{\pi\sqrt{|d|}}\zeta(2s)\left(\frac{1}{s-1} -\frac{c(s)}{s}\left(\frac{ |d|}{4}\right)^{\!1/2 -s}\right)+O(h(d)|d|^{-\delta-\sigma/2}).
\end{equation}

Suppose now that for some sequence of fundamental discriminants $d$, we can choose $\delta=\frac{1}{4}+\vep$ with arbitrary $\vep>0$  in the Linnik asymptotics  (see the corresponding arguments in \cite{VT-1}). Let $\zeta(\rho)=0$, where $\re\rho>1/2$.  Since $\zeta_{K}(s)=\zeta(s)L(s,\chi_{d})$, by letting $d\to-\infty$ we get from  \eqref{zeta-formula} that $\zeta(2\rho)=0$. This contradiction `proves' the Riemann hypothesis, as in the Deuring--Mordell--Heilbronn approach.

Of course, the starting formula \eqref{zeta-int}  needs to be proved, since for domains $\Omega$  lying on the very `top'  of the truncated modular figure $\tilde{F}_{d}$, the Linnik asymptotic expression  \eqref{uniform} loses its meaning, and the Heegner points are no longer uniformly distributed as $d\to-\infty$. The easiest way to see it is to average the representation 
 \eqref{Hecke} of a zeta function $\zeta_{D}(s)$ of order $\frak{O}_{D}$ with the discriminant $D=df^{2}$ in the imaginary quadratic field $\QQ(\sqrt{d})$ over all  
$-D\leq X$. The corresponding formula, an analogue of the classical Vinogradov-Gauss formula \cite{IMV} in the critical strip, was obtained in \cite{VT-2} and has the form
\begin{equation}\label{V-G}
\sum_{-D\leq X}\left(\frac{|D|}{4}\right)^{\!s/2}\zeta_{D}(s)=\zeta(2s)(\Phi(s)X^{1+\frac{s}{2}}+c(s)\Phi(1-s)X^{1+\frac{1-s}{2}})+R_{s}(X).
\end{equation}
Here
$$\Phi(s)=\frac{2^{-s}\zeta(s)}{(s+2)\zeta(s+2)}\quad\text{and}\quad R_{s}(X)=O\left(\frac{X^{3/4}\log^{3}X}{|s-1|(|s-1/2|+\log^{-1}X)}\right),$$
and the estimate of the remainder is uniform with respect to $s$ on compact subsets of the critical strip.
The leading term of the asymptotics in \eqref{V-G} is proportional to $\zeta(s)$, which makes the previous argument inapplicable.  Thus if Linnik's asymptotic expression can hold up to the very top
of the truncated modular figure $\tilde{F}_{d}$, then it is only for very special values of $d$. 

In the above arguments  the formula \eqref{Hecke} played the key role.  It was also used by Zagier \cite{Zagier} to construct a nontrivial representation of $\mathrm{SL}(2,\RR)$ connected with the zeros of the Riemann zeta function. We should also mention recent works by Bombieri and Garrett on the spectrum of the Laplace operator on the space of pseudo-cusp forms in connection with the zeros of $\zeta_{K}(s)$ (see the talks \cite{Bom, Gar}).  We leave the reader alone with this intriguing works and literature cited there.

\subsection*{Note added in Proof} The results of the talks [3] and [18] have now appeared on the arXiv: Enrico Bombieri and Paul Garrett, “Designed Pseudo-Laplacians”, \url{https://arxiv.org/abs/2002.07929}.

\noindent
{\bf Leon A.~Takhtajan}\\
L. Euler International Mathematical Institute;\\
Stony Brook University, USA\\
{\it E-mail}: leontak@math.stonybrook.edu


\begin{thebibliography}{1}
\bibitem{AV06} M. Aganagic, R. Dijkgraaf, A. Klemm, M. Mari\~{n}o, and C. Vafa, ``Topological strings and integrable hierarchies'', Commun. Math. Phys., \textbf{261} (2006), 451--516.
\bibitem{AG} Н.И. Ахиезер и И.М. Глазман,\emph{Теория линейных операторов в гильбертовом пространстве}, 2-е изд., Наука, М., 1966, 543 c.; English transl. N.I. Akhiezer and I.M. Glazman, \emph{Theory of linear operators in Hilbert space}, Dover Publications Inc., New York, 1993.
\bibitem{Bom} E. Bombieri, \emph{Pseudo-Laplacians: a special case}, talk at the conference ``Perspectives 
on the \\ Riemann hypothesis'' (Heilbronn Inst., 2018), \url{https://www.bristolmathsresearch.org/wp-content/uploads/2017/10/Bombieri-talk-i.pdf}.
\bibitem{BS} З.И. Боревич, И.Р. Шафаревич, \emph{Теория чисел}, 3-е доп. изд., Наука, М., 1985, 504 с.; English transl. Z.I. Borevich, I.R. Shafarevich, \emph{Number theory}, Pure and Applied Mathematics, Vol. 20 Academic Press, New York-London 1966, x+435 pp.
\bibitem{Bus-Fad} В.С. Буслаев, Л.Д. Фаддеев, “О формулах следов для дифференциального сингулярного оператора Штурма–Лиувилля”, Докл. АН СССР, \textbf{132}:1 (1960), 13--16; English transl., V.S. Buslaev, L.D. Faddeev, “Formulas for traces for a singular Sturm–Liouville differential operator”, Soviet Math. Dokl. 1 (1960), 451--454.
\bibitem{Cartier}  P. Cartier, ``Comment l'hypot\`{e}se de Riemann ne fut pas prouv\'{e}e (extraits de deux lettres de P. Cartier \`{a} A. Weil, dat\'{e}es du 12 ao\^{u}t et du 15 septembre 1979)'' (French) Seminar on Number Theory, Paris 1980-81 (Paris, 1980/1981), Progr. Math., \textbf{22} (1982), 35--48. 
 \bibitem{CV1} Y. Colin de Verdi\`{e}re, ``Pseudo-laplaciens. I'', Ann. Inst. Fourier (Grenoble) \textbf{32}:3 (1982), 275--286. 
 \bibitem{CV2} Y. Colin de Verdi\`{e}re, ``Pseudo-laplaciens. II'', Ann. Inst. Fourier (Grenoble) \textbf{33}:2 (1983), 87--113.
  \bibitem{Deuring} Max Deuring, ``Imagin\"{a}re quadratische Zahlk\"{o}rper mit der Klassenzahl 1'',  Math. Z., \textbf{37} (1933), 405--414.
 \bibitem{Dik} Л.А. Дикий, “Формулы следов для дифференциальных операторов Штурма–Лиувилля”, УМН, \textbf{13}:3(81) (1958), 111--143; English transl. L.A. Dikii, ``Trace formulas for Sturm-Liouville differential operators'', Translations AMS Series 2 \textbf{18}, (1958), 81--115.
 \bibitem{Duke} W. Duke, ``Hyperbolic distribution problems and half-integral weight Maass forms'', Invent. math., \textbf{92} (1988), 73--90.
 \bibitem{Fad59} Л.Д. Фаддеев, “Обратная задача квантовой теории рассеяния”, УМН, \textbf{14}:4(88) (1959), 57--119; English transl., L.D. Faddeev, “The inverse problem in the quantum theory of scattering”,
J. Math. Phys. 4 (1963), 72--104.
\bibitem{Fad64} Л.Д. Фаддеев, “Свойства $S$-матрицы одномерного уравнения Шредингера”, Краевые задачи математической физики. 2, Сборник работ. Посвящается памяти Владимира Андреевича Стеклова в связи со столетием со дня его рождения, Тр. МИАН СССР, \textbf{73}, Наука, М.--Л., 1964, 314--336; Amer. Math. Soc. Transl. Ser. 2, \textbf{65} (1967), 139--166; English transl., L.D. Faddeev, “Properties of the $S$-matrix of the one-dimensional Schr\"{o}dinger equation”, Amer. Math. Soc. Transl. Ser. 2, vol. 65, Amer. Math. Soc., Providence, RI 1967, pp. 139--166.
\bibitem{Fad67} Л.Д. Фаддеев, “Разложение по собственным функциям оператора Лапласа на фундаментальной области дискретной группы на плоскости Лобачевского”, Тр. ММО, \textbf{17}, Изд-во Моск. ун-та, М., 1967, 323--350; English transl., L.D. Faddeev, “Expansion in eigenfunctions of the Laplace operator on the fundamental domain of a discrete group on the Lobachevskij plane”, Trans. Moscow Math. Soc. 17 (1967), 357--386.
\bibitem{Fad74}  Л.Д. Фаддеев, “Обратная задача квантовой теории рассеяния. II”, Итоги науки и техн. Сер. Соврем. пробл. мат., \textbf{3}, ВИНИТИ, М., 1974, 93--180; English transl., L.D. Faddeev, “Inverse problem of quantum scattering theory. II”, J. Soviet Math. 5:3 (1976), 334--396.
\bibitem{F00}L.D. Faddeev, ``Modular double of a quantum group'', \emph{Conf\'{e}rence Mosh\'{e} Flato 1999, Vol. I (Dijon)}, 149--156, \emph{Math. Phys. Stud.}, \textbf{21}, Kluwer, 2000.
\bibitem{Fay} John D. Fay, ``Fourier coefficients of the resolvent for a Fuchsian group'', J. Reine Angew. Math., \textbf{293(294)} (1977), 143-203.
\bibitem{Gar} P. Garrett, \emph{Self-adjoint operators on automorphic forms}, talk at the conference ``Perspectives on the Riemann Hypothesis'' (Heilbronn Inst., 2018), \url{https://www.bristolmathsresearch.org/wp-content/uploads/2017/10/Garrett-talk-i.pdf}.
\bibitem{GD} И.М. Гельфанд, Л.А. Дикий, “Асимптотика резольвенты штурм–лиувиллевских уравнений и алгебра уравнений Кортевега–де Фриза”, УМН, \textbf{30}:5(185) (1975), 67--100;  English transl. I. M. Gel'fand, L. A. Dikii, ``Asymptotic benaviour of the resolvent of Sturm–Liouville equations and the algebra of the Korteweg--de Vries equations'', Russian Math. Surveys, 30:5 (1975), 77--113.
\bibitem{GPS} И.М. Гельфанд, М.И. Граев, И.И. Пятецкий-Шапиро, \emph{Теория представлений и автоморфные функции}, Обобщенные функции, \textbf{6}, Наука, М., 1966, 512 c.; English transl. I.M. Gel’fand, M.I. Graev, I.I. Pyatetskii-Shapiro, \emph{Representation theory and automorphic functions}, W.B. Saunders Co., Philadelphia, PA--London--Toronto, ON, 1969, xvi+426 с.
\bibitem{GL} И.М. Гельфанд и Б.М. Левитан, ``Об одном простом тождестве для собственных значений дифференциального оператора второго порядка'', ДАН CCCР \textbf{88}:4, (1953), 593—596; I.M. Gelfand, B.M. Levitan, ``On a simple identity for the eigenvalues of a second-order differential operator'', Dokl. Akad. Nauk. USSR, \textbf{88}, (1953), 953--956.
\bibitem{GK} И.Ц. Гохберг и М.Г. Крейн, \emph{Введение в теорию линейных несамосопряженных операторов в гильбертовом пространстве}, Наука, М., 1965, 448 с.; English transl., I.C. Gohberg and M.G. Krein, \emph{Introduction to the theory of linear nonselfadjoint operators}, Transl. Math. Monogr., vol. 18, Amer. Math. Soc., Providence, RI 1969, xv+378 pp.
\bibitem{GHM16} A. Grassi, Y. Hatsuda, M. Mari\~{n}o, ``Topological strings from quantum mechanics'', Annales Henri Poincar\'{e},  \textbf{17}:11 (2016), 3177--3235.
\bibitem{H} H. Heilbronn, ``On the class number in imaginary quadratic fields'', Quart. J. Math. (Oxford), \textbf{5} (1934), 150--160.
\bibitem{Hejhal0} Dennis A. Hejhal, ``Some observations concerning eigenvalues of the Laplacian and Dirichlet $L$-series'', Recent Progress in Analytic Number Theory, vol. 2, Academic Press, 1981, 95--110.
\bibitem{Hejhal} Dennis A. Hejhal, \emph{The Selberg trace formula for $\mathrm{PSL}(2,\RR)$ Vol. 2}, Lect. Notes in Math, \textbf{1001}, Springer-Verlag, Berlin, 1983.
\bibitem{Ivan} H. Iwaniec, ``Fourier coefficients of modular forms of half-integral weight'', Invent. Math,. \textbf{87} (1987), 385--401.
\bibitem{K} Felix Klein, \emph{Vorlesungen \"{ub}er die Entwicklung der Mathematik in 19.Jahrhundert}, Reprint, Springer-Verlag, Berlin Heidelberg New York, 1979.
\bibitem{Lang} S. Lang, $\mathrm{SL}_{2}(\RR)$, Addison-Wesley Publishing Co., Reading, MA-London-Amsterdam, 1975,
xvi+428 pp.
\bibitem{LST16} Ari Laptev,  Lucas Schimmer and Leon A. Takhtajan, ``Weyl type asymptotics and bounds for the eigenvalues of functional-difference operators for mirror curves'', Geom. Funct. Anal. (GAFA), \textbf{26} (2016), 288--305.
\bibitem{Lax} Peter D. Lax, \emph{Functional analysis}, Wiley-Interscience, 2002.
\bibitem{LP} P.D. Lax, R.S. Phillips, \emph{Scattering theory for automorphic functions}, Ann. of Math. Stud., vol. 87, Princeton Univ. Press, Princeton, 1976, x+300 с.
\bibitem{LS} Б.М. Левитан и И.С. Саргсян,  \emph{Введение в спектральную теорию: самосопряженные обыкновенные дифференциальные операторы}, Наука, М., 1970; English transl. B.M. Levitan, I.S. Sargsjan, \emph{Introduction to spectral theory: selfadjoint ordinary differential operators}, Translations of mathematical monographs, vol. 39, Amer. Math. Soc., Providence, RI, 1975.
\bibitem{Linnik} Ю.В. Линник, ``Асимптотическое распределение приведенных бинарных квадратичных форм в связи с геометрией Лобачевского'', Избранные труды. Теория чисел. Эргодический метод и $L$-фунцкии, Наука, Л., 1979, 141--200; Yu.V. Linnik, ``Asymptotic distribution of reduced binary quadratic forms in connection with Lobachevsky geometry'', Selected Works. Number theory. Ergodic method and $L$-functions, Nauka, L., 1979, 141--200.
\bibitem{Mar77} В.А. Марченко, \emph{Операторы Штурма-Лиувилля и их приложения}, Наукова Думка, Киев, 1977, 331 с.; English transl., V.A. Marchenko, Sturm--Liouville operators and applications, Oper. Theory Adv. Appl., vol. 22, Birkh\"{a}user Verlag, Basel 1986, xii+367 pp.
\bibitem{Mordell} L.J. Mordell, ``On the Riemann hypothesis and imaginary quadratic fields with a given class number'', J. London Math. Soc., \textbf{9} (1934), 289--298.
\bibitem{Sad}В.А. Садовничий, В.Е. Подольский, “Следы операторов”, УМН, \textbf{61}:5(371) (2006), 89–156; English transl. V.A. Sadovnichii, V.E. Podolskii, ``Traces of operators'', Russian Math. Surveys, 61:5 (2006), 885--953.
\bibitem{Selberg} A. Selberg, ``Harmonic analysis and discontinuous groups in weakly symmetric Riemannian spaces with applications to Dirichlet series'', Indian Journ. Math. Soc. \textbf{20} (1956), 47--87.
\bibitem{Schwarz} А.С. Шварц, “Эллиптические операторы в квантовой теории поля”, Итоги науки и техн. Сер. Соврем. пробл. мат., \textbf{17}, ВИНИТИ, М., 1981, 113–173; English transl. A.S. Schwarz, ``Elliptic operators in quantum field theory'', J. Soviet Math., 21:4 (1983), 551--601.
\bibitem{Siegel} Carl Ludwig Siegel, ``\"{U}ber die Classenzahl quadratischer Zahlk\"{o}rper'', Acta Arith., \textbf{1}:1 (1935), 83-86.
\bibitem{T} L.A. Takhtajan, Quantum mechanics for mathematicians, Grad. Stud. Math., \textbf{95}, Amer. Math. Soc., Providence, RI, 2008, xvi+387 pp; Russian transl., RCD, Moscow–Izhevsk 2011.
\bibitem{T17} Л.А. Тахтаджян, А.Ю. Алексеев, И.Я. Арефьева, М.А. Семенов-Тян-Шанский, Е.К. Склянин, Ф.А. Смирнов, С.Л. Шаташвили, ``Научное наследие Л.Д. Фаддеева. Обзор работ'', УМН, \textbf{72}:6(438) (2017), 3--112; English transl.  L.A. Takhtajan, A.Yu. Alekseev, I.Ya. Aref’eva, M.A. Semenov-Tian-Shansky, E.K. Sklyanin, F.A. Smirnov, and S.L. Shatashvili, ``Scientific heritage of L.D. Faddeev. Survey of papers'', Russian Math. Surveys, \textbf{72}:6 (2017), 977--1081.
\bibitem{TF15} Л.А.~Тахтаджян, Л.Д. Фаддеев, “Спектральная теория одного функционально-разно\-ст\-ного оператора конформной теории поля”, Изв. РАН. Сер. матем., \textbf{79}:2 (2015), 181--204; English transl., L.A. Takhtajan and L.D. Faddeev, “The spectral theory of a functional-difference operator in conformal field theory”, Izv. Math. 79:2 (2015), 388--410.
\bibitem{Venkov} А. Б. Венков, “Спектральная теория автоморфных функций”, Тр. МИАН СССР, \textbf{153}, 1981, 3--171; English transl. A.B. Venkov, ``Spectral theory of automorphic functions''. Proc. Steklov Inst. Math., \textbf{153} (1982), 1--163. 
\bibitem{VKF73} А.Б. Венков, В.Л. Калинин, Л.Д. Фаддеев, “Неарифметический вывод формулы следа Сельберга”, Дифференциальная геометрия, группы Ли и механика, Зап. научн. сем. ЛОМИ, \textbf{37}, Изд-во «Наука», Ленингр. отд., Л., 1973, 5--42; English transl., A.B. Venkov, V.L. Kalinin and L.D. Faddeev, “A non-arithmetic derivation of the Selberg trace formula”, J. Soviet Math. 8:2 (1977), 171–199.
\bibitem{IMV} И.М. Виноградов, ``О среднем значении числа классов чисто коренных форм отрицательного определителя'', Сообщ. Харьк. мат. об-ва, \textbf{16} (1917), 10--38. См.: Избр.
труды. М.: Изд-во АН СССР, 1952, 29--53; I.M. Vinogradov, ``On the mean value of the number of classes of properly primitive forms of negative discriminant'', Soobshch. Khark. Mat. Obshch., 16, 10--38 (1917) (Selected Works, 29--53).
\bibitem{VT-1} А. И. Виноградов, Л. А. Тахтаджян, “Об асимптотиках Линника–Скубенко”, Докл. АН СССР, \textbf{253}:4 (1980), 777--780; English transl. A.I. Vinogradov and L.A. Tahtad\v{z}jan, ``On the Linnik-Skubenko asymptotics'', Soviet  Math. Dokl. \textbf{22}:1 (1980), 136--140.
\bibitem{VT-2} А. И. Виноградов, Л. А. Тахтаджян, “Аналоги формулы Виноградова--Гаусса в критической полосе”, Аналитическая теория чисел, математический анализ и их приложения, Сборник статей. Посвящается академику Ивану Матвеевичу Виноградову к его к его девяностолетию, Тр. МИАН СССР, \textbf{158}, 1981, 45--68; English transl. A. I. Vinogradov, L. A. Takhtadzhyan, ``Analogues of the Vinogradov--Gauss formula in the critical strip'', Proc. Steklov Inst. Math., 158 (1983), 47--71.
\bibitem{Zagier} D. Zagier, ``Eisenstein series and the Riemann zeta function'', in 1979 Bombay Colloquium on Automorphic Forms, Representation Theory, and Arithmetic, Tata Inst. Fund. Research, 1981,
275--301.
\bibitem{ZF} В.Е. Захаров, Л.Д. Фаддеев, “Уравнение Кортевега--де Фриса --- вполне интегрируемая гамильтонова система”, Функц. анализ и его прил., \textbf{5}:4 (1971), 18--27; English transl., V.E. Zakharov and L.D. Faddeev, “Korteweg--de Vries equation: a completely integrable Hamiltonian system”, Funct. Anal. Appl. 5:4 (1971), 280--287.
\end{thebibliography}
\end{document}